\documentclass[a4paper,oneside,11pt]{article}
\usepackage{amsmath}
\usepackage{amssymb}
\usepackage{titlesec}
\usepackage{graphicx}
\usepackage{fancyhdr}
\usepackage{latexsym}
\usepackage{mathrsfs}
\usepackage{booktabs}
\usepackage{mathrsfs}
\usepackage{mathtools}
\usepackage{amsthm}
\usepackage{wasysym}
\usepackage{cite}
\usepackage{color}
\usepackage{amsfonts}
\usepackage{pgf,pgfarrows,pgfnodes,pgfautomata,pgfheaps}
\usepackage{amsmath,amssymb}
\usepackage{graphics}
\usepackage{multimedia}
\usepackage{graphics}
\usepackage{cite}
\usepackage{latexsym}
\usepackage{amsmath}
\usepackage{amssymb}
\usepackage[dvips]{epsfig}
\usepackage{amscd}
\numberwithin{equation}{section}

\newtheorem{theorem}{Theorem}[section]
\newtheorem{lemma}{Lemma}[section]
\newtheorem{remark}{Remark}[section]

\newcommand{\thatsall}{\hfill$\Box$}
\newcommand{\R}{{\mathbb{R}^2}}

\date{}
\title{Global Well-posedness of Incompressible MHD System with Density-dependent Viscosity in $\mathbb{R}^2$ }
\author{ Xue W{\small ANG}  \thanks{Email:  xuewa@amss.ac.cn (X.  Wang) } \\
{\normalsize  School of Mathematical  Sciences,} \\ {\normalsize    Tiangong University, Tianjin 300387, P. R. China}
}

\begin{document}
\maketitle
\begin{abstract}

We study the Cauchy problem for the two-dimensional
inhomogeneous incompressible Magnetohydrodynamic equations
with density-dependent viscosity and vacuum. For initial
data with positive mass, finite second moment of density, and zero total momentum, we establish the existence
and uniqueness of global strong solutions when the viscosity is sufficiently close to one in $L^\infty$. The smallness threshold may depend on the initial data, while no separate smallness condition is imposed on the initial velocity or magnetic field. We obtain uniform-in-time bounds for the density second moment and its $H^1\cap W^{1,q}$ norm, where $q>2$. Moreover, for every $\eta\in(0,1)$, $\|\sqrt{\rho}u(\cdot,t)\|_{L^2}+\|b(\cdot,t)\|_{L^2}$ decays at the rate $t^{-\eta/2}$, and the $H^1$ norms of
$\nabla u$, $P$, and $\nabla b$ decay at the rate
$t^{-(1+\eta)/2}$. The proof combines a velocity decomposition based on Hardy-space estimates with magnetic-field integrability and time-weighted energy estimates.

\end{abstract}

\section{Introduction}
\quad We consider an inhomogeneous incompressible viscous fluid in the whole space $\mathbb{R}^2$ governed by the  Magnetohydrodynamic (MHD) equations: 
\begin{equation}\label{C101}
\begin{cases}
\rho_t+\mathrm{div}(\rho u)=0,\\
(\rho u)_t+\mathrm{div}(\rho u \otimes u)-\mathrm{div}(2\mu(\rho)Du)+\nabla P-b\cdot\nabla b=0,\\
b_t+u\cdot\nabla b-b\cdot\nabla u-\nu\Delta b=0,\\
\mathrm{div}u=0, \quad \mathrm{div}b=0,
\end{cases}
\end{equation}
where $t\ge 0$, $x=(x^1,x^2)\in\mathbb{R}^2$, the unknowns  $\rho$, $u=(u^1,u^2)$, $P$, and $b=(b^1,b^2)$ stand for the fluid density, velocity field, pressure and magnetic field, respectively,  $Du=( \nabla u+(\nabla u)^T)/2$ is the deformation tensor. In general, the viscosity $\mu =\mu(\rho)$ is a function of $\rho$ satisfying 
\begin{equation} \label{mu}
\mu(s) \in C^1[0, \infty) ~~ \mbox{and} ~~  \mu(s)>0, ~~ \text{for  any }\,s\in[0,\infty).
\end{equation}
The resistivity $\nu$ is a positive constant. Without loss of generality, we assume that $\nu=1$ in the following paper. 

To study the Cauchy problem for the equations  \eqref{C101} with $(\rho,b)$ vanishing at infinity in some weak sense, we impose \eqref{C101} the initial conditions:
\begin{equation}\label{1.2}
\rho(x,0)=\rho_0(x),\quad \rho u(x,0)=\rho_0 u_0(x),\quad b(x,0)=b_0(x),\quad x\in\R.
\end{equation}

When the magnetic field vanishes, i.e., $b\equiv 0$,  equations \eqref{C101} reduce to the nonhomogeneous Navier-Stokes equations with variable viscosity. Since extensive results are already available for that case (see, e.g., \cite{LS,K,Sm,DM2,hll} and the references therein), we shall not review them in detail here and instead focus on the MHD system.

In the homogeneous case, where $\rho$ is a constant, the MHD system has been studied widely. Duvaut-Lions  \cite{DL} established the global existence of weak solutions with finite energy and local strong solutions; in particular, the 2D local strong solution was proved to be global and unique, see also \cite{ST}. For the three-dimensional case, various Serrin-type criteria analogous to those for the Navier-Stokes equations were established in \cite{CMZ,HX,wjh}. As for the 3D Navier-Stokes equations, whether the local strong solution is global remains an outstanding open problem for the homogeneous MHD system as well.

When the fluid is nonhomogeneous, the situation becomes more complicated. The finite-energy weak solutions in 3D bounded domains is constructed by Gerbeau-Le Bris \cite{GB}, the corresponding result in the torus is given by Desjardins-Le Bris \cite{DB}. In the frame of Besov spaces, Abidi et al.  \cite{AH,AP} obtained the global well-posedness of strong solutions to the 3D Cauchy problem under the assumptions that the initial velocity and magnetic field are sufficiently small and that the initial density is close to a positive constant, so that no vacuum is present. Chen-Tan-Wang \cite{CTW} extended the local existence result in \cite{AH}  to the case where vacuum is allowed. In addition, in the absence of vacuum, Sokrani \cite{So} proved the global existence of strong solutions for small initial data in certain Sobolev spaces, generalizing the corresponding result for nonhomogeneous Navier–Stokes equations with variable viscosity obtained by Abidi-Zhang \cite{AZ}.

When vacuum is present, the problem becomes even more delicate. Under a natural compatibility condition, 
Li \cite{Li} established the global-in-time unique strong solution with density-dependent viscosity and resistivity  for the 3D case, provided that $\|\nabla u_0\|_{L^2}+\|\nabla b_0\|_{L^2}$ is sufficiently small. This result was later improved by Liu \cite{Liu}, who removed the compatibility condition under the smallness assumption $\|\rho_0\|_{L^\infty}+\|b_0\|_{L^3}$ suitably small. Very recently, Zhang \cite{Z} investigated the global existence and large-time asymptotic behavior of strong solutions to the 3D Cauchy problem, provided that the initial velocity and magnetic field are suitably small in the $\dot{H}^\beta$-norm for some $\beta\in(\frac12,1]$. 
For the 3D case in exterior domains, Yuan-Zhang-Zhou \cite{YZZ0} obtained the global strong solutions with density-dependent viscosity under slip boundary conditions. For the general large data, Huang-Wang \cite{HW} and L\"u-Xu-Zhong \cite{LXZ} showed the global existence of strong solutions with constant viscosity to the initial boundary value problem in 2D bounded domains and the Cauchy problem  in $\mathbb{R}^2$, respectively. In the 2D setting with density-dependent viscosity, Zhong \cite{zx} established the global well-posedness and  exponential decay of strong solutions in bounded domains  provided that $\|\nabla\mu(\rho_0)\|_{L^q}$ $(q>2)$ is suitably small.  

As pointed out in \cite{hll}, the strong interaction between the density and the velocity field introduces serious difficulties in the mathematical study of the global theory for the case of density-dependent viscosity, and the methods developed for the constant-viscosity case cannot be applied directly. Furthermore, the strong coupling between the fluid motion and the magnetic field brings further challenges, there needs some new ideas compared to the case of Navier-Stokes provided by \cite{FLW}.  In the present paper, we aim to investigate the global existence and  decay rates  of strong solutions to the problem \eqref{C101}--\eqref{1.2}, via time-weighted techniques and energy methods.

Before stating the main results, we first explain the notations and conventions used throughout the paper. We denote the spaces of Schwartz functions and tempered distributions in $\R$ by $\mathcal{S}$ and $\mathcal{S}'$,  respectively. For the integer $k\geq 0$ and $1 \leq p \le \infty$, the  standard Lebesgue spaces, H\"{o}lder spaces, and Sobolev ones are denoted by:
\begin{equation*}\begin{cases}
L^p\triangleq L^p(\mathbb{R}^2),\quad C^\alpha\triangleq C^\alpha(\R),\quad  W^{k, p}\triangleq W^{k, p}(\mathbb{R}^2),\quad H^k\triangleq W^{k, 2},\\
D^{k,p}\triangleq D^{k,p}(\mathbb{R}^2)=\{v\in L^1_{\mathrm{loc}}(\mathbb{R}^2)\,\big|\,\nabla^k v\in L^p\},\quad D^k\triangleq D^{k,2},\\
C_{0,0}^\infty=\big\{w\in C_0^\infty(\mathbb{R}^2)\,\big|\,\mathrm{div}w=0\big\},\quad D^1_{0,0}= \mathrm{the~closure~of}~{C_{0,0}^\infty}~\mathrm{in}~D^1.
\end{cases}
\end{equation*}
Set $\tilde x\triangleq (1+|x|^2)^{1/2}$ and   $B_N\triangleq \{x\in \mathbb{R}^2\big|\,|x|<N\}$ for $N>0$. The material derivative of a function $f$ is represented by $\dot{f}\triangleq\partial_tf+u\cdot\nabla f.$




\begin{theorem}\label{T1}
For given constants $\bar\rho>0$ and $q>2$, suppose that the initial data $(\rho_0, u_0,b_0)$ satisfy
\begin{equation}\label{1.3}
\begin{split}
 0&\leq \rho_0\leq \bar\rho,\quad \rho_0\,\tilde x^2 \in L^1,\quad  \rho_0\in H^1\cap W^{1, q},\\
 \sqrt{\rho_0}u_0&\in L^2,\quad u_0\in D^1_{0,0},\quad b_0\in L^1\cap H^1,\quad \mathrm{div}b_0=0,\\
 0<&M\triangleq\int_{\mathbb{R}^2}\rho_0(x)\,dx<\infty,\quad\int_{\mathbb{R}^2}\rho_0u_0(x)\,dx=0.
\end{split}
\end{equation}
Then there exists a positive constant $\varepsilon_0$ depending on $\bar\rho$, $\mu$, $q$, and the initial data 
such that if 
\begin{equation}\label{com}
\|\mu(\rho_0)-1\|_{L^\infty}\le \varepsilon_0,
\end{equation}
the system \eqref{C101}--\eqref{1.2} admits a unique global strong solution $(\rho, u,P,b)$ satisfying that for any  $0<\tau<T<\infty$ and $1<p<\infty$,
\begin{equation}\begin{cases}\label{1.7} 
\rho\in  C([0,T];L^1\cap H^1\cap W^{1, q}),\quad
\rho\,\tilde x^2\in L^\infty(0,T;L^1),\\
\sqrt\rho u,\nabla u,\nabla b\in L^\infty(0,T;L^2),\quad b\in L^\infty(0,T;L^p);\\
\sqrt\rho u_t\in L^2(0,T;L^2)\cap L^\infty(\tau,T;L^2),\quad \nabla u_t\in L^2(\tau,T;L^2),\\
\nabla u,P\in L^2(0,T;H^1)\cap L^\infty(\tau,T;H^1)\cap L^2(\tau,T;W^{1,q}).\\
b\in L^\infty(\tau,T;H^2)\cap L^2(\tau,T;H^{3}),\quad b_t\in L^\infty(\tau,T;L^2)\cap L^2(\tau,T;H^{1}).
\end{cases}\end{equation}

Moreover, the strong solution $(\rho,u,P,b)$ has the following large time behaviors:

1. $($Uniform bounds on $\rho)$  There is a positive constant $\mathbf{C}_1$ decided by $\bar\rho$, $\mu$, $q$, and the initial values 
such that for any $T>0$,
\begin{equation}\label{C102}
\sup_{0\leq t\leq T}\|\rho(\cdot,t)\|_{L^\infty}\leq\bar\rho,\quad \sup_{0\leq t\leq T}\left(\|\rho(\cdot,t)\,\tilde x^2\|_{L^1}+\|\rho(\cdot,t)\|_{ H^1\cap W^{1,q}}\right)\leq\mathbf{C}_1.
\end{equation}

2. $($Algebraic decay of $u$, $P$, and $b)$ For any $\eta\in(0,1)$, there is a positive constant $\mathbf{C}_2$ depending on $\eta$, $\bar\rho$, $\mu$, $q$, and the initial values 
such that  
for any $t\geq 1$, 
\begin{equation}\begin{aligned}\label{1.8}
&\|\sqrt\rho u(\cdot,t)\|_{L^2}+\|b(\cdot,t)\|_{L^2}\le \mathbf{C}_2\,t^{-\eta/2},\\
 \|\nabla u(\cdot,t)&\|_{H^1}+\|P(\cdot,t)\|_{H^{1}}+\|\nabla b(\cdot,t)\|_{H^1}
\leq \mathbf{C}_2\,t^{-(1+\eta)/2}.
\end{aligned}\end{equation}

3. (Far field state of $u$) For a.e. $t>0$, there exists a unique far field velocity $u_\infty(t)$ such that $u(x,t)\rightarrow u_\infty(t)$ as $|x|\rightarrow\infty$, and for any  $\alpha\in(0,1/2)$ and $\eta\in(0,1)$, there is a positive  constant $\mathbf C_3$  depending on $\alpha$, $\eta$, $\bar\rho$, $\mu$,  and the initial values such that 
\begin{equation}\label{1.9}
\| u-u_\infty\|_{L^{2/\alpha}}+|u_\infty(t)|\le \mathbf{ C}_3\,t^{-(1+\eta)/2}.
\end{equation}
\end{theorem}
\begin{remark}
The assumptions in Theorem \ref{T1} have distinct roles.
The positive mass allows the spatially constant component
of the velocity to be determined from the total momentum.
The zero-momentum condition selects the reference frame
in which the kinetic energy decays to zero.
The finite second moment of the density is used in the
Hardy-space and weighted elliptic estimates, while
$b_0\in L^1\cap H^1$ provides the initial $L^p$ bounds,
$1<p<2$, used to derive the magnetic-field decay. For non-zero initial total momentum, we can use the Galilean transformation reducing the momentum to zero


The decay estimates are established for every
$\eta\in(0,1)$. The endpoint $\eta=1$ and the optimality of
these rates are not addressed here.
\end{remark}



We now explain the main difficulties and the ideas of the proof. Our analysis builds on the Hardy-space approach of
Fan-L\"u-Wang \cite{FLW}. The magnetic coupling introduces an additional difficulty in the large-time analysis. Indeed, the velocity equation contains the magnetic stress
$\mathrm{div}(b\otimes b)$. We therefore write
$$u=w_1+w_2+u_\infty,$$
where $w_1$ and $w_2$ correspond to the inertial term
$\rho\dot u$ and the magnetic stress, respectively. Here, $u_\infty$ is a spatially constant vector. The cancellation property $\int_{\mathbb R^2}\rho\dot u\,dx=0$
allows us to estimate $w_1$ by Hardy-space and weighted
elliptic estimates. The component $w_2$ is controlled by the integrability of the magnetic field. The zero-momentum condition determines $u_\infty$.

Set
$$M_2(t)
\triangleq\int_{\mathbb R^2}\rho(x,t)\,(1+|x|^2)\,dx.
$$
The decomposition yields
$$
\|\sqrt{\rho}\,u\|_{L^2}^2\leq C\,(1+M_2(t))^{1/3}\|\sqrt{\rho}\dot u\|_{L^2}^2+C\,\|\nabla b\|_{L^2}^2.
$$
In contrast to the Navier-Stokes case, this estimate
retains a magnetic contribution.
The total energy also contains $\|b\|_{L^2}^2$.
On $\mathbb R^2$, this quantity cannot be controlled by
$\|\nabla b\|_{L^2}^2$ through a Poincar\'e inequality.
Thus, the energy closure used in the Navier-Stokes
exponential-decay argument does not carry over directly.
We instead propagate the magnetic-field $L^p$ bounds
for $1<p<2$. These bounds yield magnetic-field decay.
They also allow us to close a hierarchy of time-weighted
estimates for the coupled system.

A further difficulty concerns the density second moment.
The transport equation gives
$$
\frac{d}{dt}M_2(t)^{1/2}\leq C\|\sqrt{\rho}u(\cdot,t)\|_{L^2}.
$$
However, the algebraic energy estimate
$$
\|\sqrt{\rho}u(\cdot,t)\|_{L^2}\leq Ct^{-\eta/2},
\quad \eta\in(0,1),
$$
does not provide an integrable upper bound on $(1,\infty)$.
It therefore does not suffice to obtain a uniform bound
for $M_2(t)$.
The distinction between energy decay and time integrability
is essential at this step.

To overcome this difficulty, we retain the precise magnetic terms in the velocity estimate. More specifically, we prove
$$
\|u-u_\infty\|_{L^8}\leq
C\,(1+M_2(t))^{1/6} \|\sqrt{\rho}\dot u\|_{L^2}+C\left(
\|b\|_{L^{16/15}}^{2/3} \|\nabla b\|_{L^2}^{4/3}+\|b\|_{L^2}\|\nabla b\|_{L^2}\right).
$$
The zero-momentum condition implies
$$
\|\sqrt{\rho}u\|_{L^2}\leq\|\rho\|_{L^{4/3}}^{1/2}\|u-u_\infty\|_{L^8},
$$
which together with preliminary moment bound $M_2(t)\leq C\,(1+t)^2$,  the time-weighted estimates, and the magnetic-field decay gives
$$
\|\sqrt{\rho}u(\cdot,t)\|_{L^2}\leq Ct^{-1-\delta},
\quad t\geq1,
$$
for some $\delta>0$. Consequently,
$$
\int_0^\infty\|\sqrt{\rho}u(t)\|_{L^2}\,dt<\infty,
\quad \sup_{t\geq0}M_2(t)<\infty.
$$
Thus, the preliminary growth estimate for the density
second moment improves to a uniform-in-time bound.

The Hardy-space estimates are inherited from \cite{FLW}.
The higher-order magnetic estimates use established MHD techniques. The additional step is to connect these estimates through time-integrable bounds for the weighted velocity. The precise magnetic terms above are crucial to this step. They allow us to control the density second moment despite the magnetic coupling. We also establish
$\nabla u\in L^1(0,\infty;L^\infty)$. This yields uniform bounds for the density gradients. Together with the higher-order estimates, these bounds close the global continuation argument. The resulting conclusions combine algebraic decay with uniform control of the density second moment and its $H^1\cap W^{1,q}$ norm.

The rest of the paper is organized as follows. In the next section, some useful basic lemmas are given. Then all the necessary a priori estimates are obtained in Section 3. In
the last section, the main result is proved.

\section{Preliminary}\label{S2}
\quad  
In this section, we collect some known facts and basic inequalities, which will be used frequently in later analysis. 

To begin with, we give the result on the local existence of the strong solutions to MHD system \eqref{C101}--\eqref{1.2} without additional weighted  conditions of $\nabla \rho_0$ and $b_0$. This can be proved in standard arguments as \cite{LXZ1,Song} and with the help of the  integrability of the difference of $u$ and its far field state, which can be derived in a similar manner as this paper,  we omit its proof here.  
\begin{lemma}\label{L41}
Under the assumptions of Theorem \ref{T1}, with the requirement $b_0\in L^1$ omitted,  there exists a positive time $T$ determined by initial data, such that the system \eqref{C101}--\eqref{1.2}  admits a unique strong solution $(\rho,u,P,b)$ in $\R\times[0,T]$ satisfying the regularity \eqref{1.7} with $2\le p<\infty$. 
\end{lemma}

The following classical Caffarelli-Kohn-Nirenberg inequality (see \cite{ckn,cw}) is the key to obtaining the integrability of $u$. 
\begin{lemma}\label{pl2}
For any function $v\in C_0^\infty$ and  real numbers $\alpha$, $\beta$ satisfying
$$\beta>0,\quad \beta-1\leq\alpha<\beta,\quad 
q=\frac{2}{\beta-\alpha},$$ 
 there exists a positive  constant $C$ determined by $\alpha$ and $\beta$ such that 
\begin{equation*}
\|v\,|x|^{\alpha}\|_{L^q}\leq C\,\|\nabla v \,|x|^{\beta}\|_{L^2}.
\end{equation*}
In fact, by density arguments, it holds for any $v$ satisfying $\nabla v\,|x|^\beta\in L^2$, $v|x|^\alpha\in L^q$.
\end{lemma}

Recall the Riesz potential $I_s(f)$ defined for $f\in\mathcal{S}$ and $0<s<2$,
$$I_s(f)(x)\triangleq\int_\R\frac{1}{|x-y|^{2-s}}f(y)\,dy.$$
We show the following weighted Hardy-Littlewood-Sobolev inequality, initially introduced by Stein-Weiss \cite{SGW}.  
\begin{lemma}\label{hls} For any  $\alpha$,  $\beta$, and $1< p\leq q<\infty$ satisfying 
\begin{equation*} 
\begin{split}
q\alpha>-2,\quad
\frac{1}{q}=\frac{1}{p}+\frac{\beta-\alpha-s}{2},\quad
\alpha\leq\beta,\quad
\beta<\frac{2(p-1)}{p},
\end{split}
\end{equation*}
there is a generic positive constant $C$ determined by $p$, $s$, $\alpha$, and $\beta$ such that for $f\in\mathcal{S}$,
\begin{equation}
\label{C307}
\|I_sf\,|x|^{\alpha}\|_{L^{q}}\leq C\,\|f\,|x|^{\beta}\|_{L^p}.
\end{equation}
\end{lemma}

 Next, we have the following Poincar$\acute{e}$-type inequality for functions in $D^{1}(\R)$, which is provided by \cite [Theorem B.1]{lions}.
\begin{lemma}\label{w3}
For $m\in[2,\infty)$ and $\zeta\in(1+m/2,\infty)$, there exists a positive constant $C=C(m,\zeta)$ such that for any $v\in D^{1}(\R)$,
\begin{equation}\label{e2}
\left(\int_\R\frac{|v|^m}{e+|x|^2}(\log(e+|x|^2))^{-\zeta}dx\right)^{1/m}\le C\left(\|v\|_{L^2(B_1)}+\|\nabla v\|_{L^2}\right).
\end{equation}
\end{lemma}

Based on this, the following lemma gives a general integrability of $u$ (see also \cite{lx}).
\begin{lemma}\label{jq} 
Suppose that $\rho$ is a non-negative function in $\R$ satisfying
\begin{equation}\label{con}
\|\rho\|_{L^1\cap L^\infty}\leq M_1,\quad \int_{B_{N}}\rho\, dx\geq M_2,\quad \int_\R\rho\,\tilde x^2\,dx\leq M_3, 
\end{equation}
with $N\geq 1$ and $M_1,M_2,M_3\in(0,\infty)$. Then for any $p\in [2,\infty)$ and $\xi\in(0,3/(5p))$, there exists a positive constant $C$ determined by $M_1$, $M_2$,  $p$, and $\xi$ such that for any $v\in C_0^\infty(\R)$,
\begin{equation}\label{2.1}
\|\rho v\|_{L^p}\leq CM_3^\xi N^2\left(\|\sqrt\rho v\|_{L^2}+\|\nabla v\|_{L^2}\right).
\end{equation}
\end{lemma}
\begin{proof}
Denote $ v_{B_N}\triangleq \frac{1}{|B_N|}\int_{B_N}v(x)\,dx$, then the Poincare's inequality implies
\begin{equation}\begin{aligned}\notag
|v_{B_N}|^2\int_{B_{N}}\rho\,dx
&\le 2\int_{B_{N}}\rho |v|^2\,dx+2\int_{B_{N}}\rho |v-v_{B_N}|^2\,dx\\
&\le 2\,\|\sqrt\rho v\|_{L^2}^2+CN^2\,\|\nabla v\|_{L^2}^2,
\end{aligned}\end{equation}
which leads to
\begin{equation}\begin{aligned}\label{local}
|v_{B_N}|\le C\,\|\sqrt\rho v\|_{L^2}+CN\,\|\nabla v\|_{L^2}.
\end{aligned}\end{equation}
Moreover, using the Poincare's inequality again and \eqref{local}, we get
\begin{equation}\begin{aligned}\notag
\|v\|_{L^2(B_N)}
&\le 2\,\|v-v_{B_N}\|_{L^2(B_N)}+C\,|B_N|^{1/2}| v_{B_N}|\\
&\le CN^2\,\|\nabla v\|_{L^2}+CN\,\|\sqrt\rho v\|_{L^2}.
\end{aligned}\end{equation}
This combined with \eqref{con} and \eqref{e2} gives that for any $\xi<3/(5p)$,
\begin{equation}\begin{aligned}\notag
\|\rho v\|_{L^p(\R)}
&\le \|\rho\,\tilde x^2\|_{L^1}^{\xi}\|v\,\tilde x^{-2\xi}\|_{L^{2/\xi}}\|\rho\|_{L^{2p(1-\xi)/(2-3p\xi)}}^{1-\xi}\\
&\le CM_3^\xi\left(\|v\|_{L^2(B_N)}+\|\nabla v\|_{L^2}\right)\\
&\le CM_3^\xi N^2\left(\|\sqrt\rho v\|_{L^2}+\|\nabla v\|_{L^2}\right).
\end{aligned}\end{equation}
The proof of this lemma is completed.
\end{proof}

Moreover, we list an important property of Hardy space $\mathcal{H}^1(\R)$, see \cite[Chapter III]{S} for the detailed proof.
\begin{lemma}\label{L21}
 For $i=1,2$, the Riesz transform $R_i$ (and similar singular integral operators, see \cite[p.71]{S}) maps $\mathcal{H}^1(\R)$ to itself, that is,  there is a positive constant $C$  such that for any $f\in\mathcal{H}^1(\R)$,
\begin{equation}\label{C202}
\|R_i(f)\|_{L^1(\mathbb{R}^2)}\leq 
 C\|f\|_{\mathcal{H}^1(\R)}.
\end{equation}
\end{lemma}

Finally, in order to estimate the derivatives of velocity,  the following regularity results for the generalized Stokes system (see \cite{FLW}) is needed.
\begin{lemma}\label{pl3}
Suppose that the function $\mu$ satisfies $\underline{\mu}\leq\mu(x)\leq\overline{\mu}$ for some $\underline{\mu},\,\overline{\mu}\in(0,\infty)$ and any $x\in\R$, we  consider the Stokes-like system in $\R$ defined by
\begin{equation}\label{2.2}\begin{cases}
-\mathrm{div}(2\mu\, Du) +\nabla P=f,\\
\mathrm{div}u=0.\\
\end{cases}
\end{equation} 
If $u$, $P$, $f$ $\in C_0^\infty(\R)$ and satisfy \eqref{2.2}, we have the following declarations.

$1).$ For any $p\in[2,\infty)$, there are positive constants $\varepsilon$ and $C$ depending on 
$\underline{\mu}$, $\overline{\mu}$, and $p$ such that if $\|\mu-1\|_{L^\infty}\leq\varepsilon$, it holds that
\begin{equation}\label{C218}
\|\nabla u\|_{L^p}\leq C\,\|\nabla u\|_{L^2}^{2/p}\|f\|_{L^2}^{1-2/p}.
\end{equation}
 
$2).$
For any $p\in(1,\infty)$, there are positive constants $\varepsilon$ and $C$ depending on  $\underline\mu$, $\overline\mu$, and $p$ such that if $\|\mu-1\|_{L^\infty}\leq\varepsilon$, it holds that
\begin{equation}\label{2.5}
\|\nabla^2 u\|_{L^p}+\|\nabla P\|_{L^p}\leq C\left(\|f\|_{L^p}+\|\nabla\mu\cdot Du\|_{L^p}\right).
\end{equation}
\end{lemma}

\section{Uniform estimates}\label{S3}
\quad  In this section, based on the local existence of strong solutions to the Cauchy problem  \eqref{C101}--\eqref{1.2}  stated in Lemma \ref{L41},  we will provide all the necessary global a priori estimates of the solutions. In the following section, we assume $(\rho,u,P,b)$ is a strong solution to \eqref{C101}--\eqref{1.2} in $\R\times[0,T]$.



We start with the following uniform estimates on density, which is guaranteed by the incompressible flow (see \cite[Lemma 3.4]{lions}).
\begin{lemma}\label{le1}
Under the conditions of Theorem \ref{T1}, we have 
\begin{equation}\begin{aligned}\label{3.7} 
\|\rho(&\cdot,t)\|_{L^\infty}=\|\rho_0\|_{L^\infty}\leq\bar\rho,\quad 
\|\rho(\cdot,t)\|_{L^1}=\|\rho_0\|_{L^1}=M,\\
&\|\mu(\rho(\cdot,t))-1\|_{L^\infty}=\|\mu(\rho_0)-1\|_{L^\infty}, \quad t\in(0,T].
\end{aligned}\end{equation}
Moreover, there exist positive constants   $\underline{\mu}$ and $\overline{\mu}$ depending on  $\bar\rho$ and $\mu$  such that 
\begin{equation}\label{1.5}
\begin{split}
0<\underline{\mu}\leq\mu(\rho(x,t))\leq\overline{\mu}<\infty,\quad (x,t)\in\R\times[0,T].\\
\end{split}
\end{equation}
\end{lemma}

\begin{lemma}\label{le2} 
Under the conditions of Theorem \ref{T1}, there exists a positive constant $C$ depending on  $\bar\rho$, $\mu$, $\|\sqrt{\rho_0}u_0\|_{L^2}$, and $\|b_0\|_{L^2}$ such that
\begin{equation}\label{3.8} 
\sup_{t\in[0,T]}\left(\|\sqrt\rho u\|_{L^2}^2+\|b\|_{L^2}^2\right)+\int_0^T\left(\|\nabla u\|_{L^2}^2+\|\nabla b\|_{L^2}^2\right)\,dt\leq C.
\end{equation}
\end{lemma}
\begin{proof}
Testing \eqref{C101}$_2$ and \eqref{C101}$_3$ by $u$ and $b$ respectively, then adding them together and using \eqref{C101}$_4$ yield   
\begin{equation}\label{3.12}
\frac{1}{2}\,\frac{d}{dt}\int_\R\left(\rho|u|^2+|b|^2\right)\,dx+\int_\R\left(2\mu(\rho)|Du|^2+|\nabla b|^2\right)\,dx=0.
\end{equation}
Since $2\|Du\|_{L^2}^2=\|\nabla u\|_{L^2}^2$ due to $\mathrm{div}u=0$, we deduce from \eqref{1.5} and \eqref{3.12} that  \eqref{3.8} holds. The proof of Lemma \ref{le2} is finished.
\end{proof}

\begin{lemma}\label{le3}
Under the conditions of Theorem \ref{T1}, suppose that 
$$\|\mu(\rho_0)-1\|_{L^\infty}\leq \varepsilon_1,$$
with $\varepsilon_1>0$ depending only on  $\bar\rho$ and $\mu$, then there exists a positive constant $C$ depending on $\bar\rho$, $\mu$, $\|\sqrt{\rho_0}u_0\|_{L^2}$,  $\|\nabla u_0\|_{L^2}$, and $\|b_0\|_{H^1}$  such that for $i=0,1$,
\begin{equation}\label{q166}
\sup_{t\in[0,T]}t^i\,\|b\|_{L^4}^4+\int_0^Tt^i\,\||\nabla b|\,|b|\|_{L^2}^2\,dt\le C,
\end{equation}
\begin{equation}\label{3.9}
\sup_{t\in[0,T]}t^i\left(\|\nabla u\|_{L^2}^2+\|\nabla b\|_{L^2}^2\right)+\int_0^T t^i\left(\|\sqrt\rho\dot u\|_{L^2}^2+\|\Delta b\|_{L^2}^2\right)\,dt\leq C.
\end{equation}
\end{lemma}
\begin{proof}
First, integrating \eqref{C101}$_3$ multiplied by $4|b|^{2} b$  by parts gives  
\begin{equation}\begin{aligned}\label{p24}
&\frac{d}{dt}\int_\R|b|^4\,dx+4\int_\R|b|^{2}|\nabla b|^2\,dx+2\int_\R|\nabla|b|^{2}|^2\,dx\\
&\le 4\|\nabla u\|_{L^2}\||b|^{2}\|_{L^4}^2\\
&\le \|\nabla |b|^{2}\|_{L^2}^2+C\,\|\nabla u\|_{L^2}^2\|b\|_{L^4}^4.
\end{aligned}\end{equation}
Multiplying this inequality by $t^i$ $(i=0,1)$ and using Gronwall's inequality and \eqref{3.8}, we get \eqref{q166} directly. 

Next, taking the $L^2$ inner product of \eqref{C101}$_2$ with $\dot{u}$ and using \eqref{C101}$_1$ yield 
\begin{equation}\label{3.6}
2\int_\R\mu(\rho)Du:D\dot u\,dx+\int_\R\rho|\dot{u}|^{2}\,dx=\int_\R P\,\mathrm{div}\dot{u}\,dx+\int_\R b\cdot\nabla b\cdot\dot u\,dx.
\end{equation}

For the first term in \eqref{3.6}, it follows from  $\mathrm{div}u=0$ and \eqref{1.5} that
\begin{equation}
\begin{split}\label{3.28}
&2\int_\R\mu(\rho) Du:D\dot u\, dx\\ 
&=2\int_\R\mu(\rho)(Du:D u_t+Du: (u\cdot\nabla Du))\,dx+2\int_\R\mu(\rho)\,\partial_iu^j\,\partial_i u\cdot\nabla u^j\,dx\\
&\geq\frac{d}{dt}\int_\R\mu(\rho)|Du|^2\,dx-\int_\R|Du|^2\left(\partial_t\mu(\rho)+u\cdot\nabla\mu(\rho)\right)\,dx-C\,\|\nabla u\|_{L^3}^3\\
&=\frac{d}{dt}\int_\R\mu(\rho)|Du|^2\,dx-C\,\|\nabla u\|_{L^3}^3,
\end{split}
\end{equation}
where in the last equality we have used 
$\partial_t\mu(\rho)+u\cdot \nabla \mu(\rho)=0$
due to \eqref{C101}.  

To estimate the term concerning $P$ in \eqref{3.6}, from \eqref{C101}$_2$, we can write $P$ as
\begin{equation}\label{C225}
P=-\Delta^{-1}\mathrm{div}(\rho\dot u-b\cdot\nabla b)+\Delta^{-1}\mathrm{div}\mathrm{div}(2\mu(\rho) Du).
\end{equation}
Owing to  $\mathrm{div}u=0$, we can see that
\begin{equation}\label{p1}
\mathrm{div}\dot u=\mathrm{div}u_t+u\cdot\nabla\mathrm{div}u+\partial_{i}u\cdot\nabla u^{i}=\partial_{i}u\cdot\nabla u^{i}.
\end{equation}
Then by \eqref{C225}, \eqref{p1}, \eqref{1.5}, and the boundedness of operators  ``$\Delta^{-1}\mathrm{div}\mathrm{div}$"   and ``$\nabla\Delta^{-1}\mathrm{div}$" in $L^p$ for $1<p<\infty$ (see \cite[Page 13]{S}), we have 
\begin{equation}\label{C416}
\begin{split}
\int_\R P\,\mathrm{div}\dot{u}\,dx
&=-\int_\R\Delta^{-1}\mathrm{div}(\rho\dot u-b\cdot\nabla b)\,\partial_{i}u\cdot\nabla u^{i}\,dx\\
&\quad+\int_\R \Delta^{-1}\mathrm{div}\mathrm{div}(2\mu(\rho)Du)\,\partial_{i}u\cdot\nabla u^i\,dx\\
&\leq \|\Delta^{-1}\mathrm{div}(\rho\dot u-b\cdot\nabla b)\|_{\mathrm{BMO}}\|\partial_{i}u\cdot\nabla u^{i}\|_{\mathcal{H}^{1}}+C\,\|\nabla u\|_{L^3}^3\\
&\leq C\,\|\nabla\Delta^{-1}\mathrm{div}(\rho\dot u-b\cdot\nabla b)\|_{L^2}\|\nabla u \|_{L^2}^{2}+C\,\|\nabla u\|_{L^3}^3\\
&\leq C\,\|\rho\dot u-b\cdot\nabla b\|_{L^2}\|\nabla u\|_{L^2}^2+C\,\|\nabla u\|_{L^3}^3.
\end{split}
\end{equation}
The main ideas (see \cite{des}) in this estimate are using  the fact: BMO$(\mathbb{R}^2)$ is the dual space of $\mathcal{H}^1(\mathbb{R}^2)$ and the div-curl lemma  (\cite[Theorem II.1]{clms}): for any $f,g\in\mathcal{D}'$, if $\mathrm{div}f=0$ and $ \mathrm{curl}g=0,$ there is a positive constant $C$ such that
\begin{equation}\notag
\|fg\|_{\mathcal{H}^1}\leq C\,\|f\|_{L^2}\|g\|_{L^2}.
\end{equation}

To deal with the last term in \eqref{3.6}, the integration by parts combined with \eqref{C101}$_4$ leads to that for any $\delta\in(0,1)$,
\begin{equation}\begin{aligned}\label{p2}
&\int_\R b\cdot\nabla b\cdot\dot u\,dx\\
&=\int_\R b\cdot\nabla b\cdot u_t \,dx+\int_\R b\cdot\nabla b\cdot (u\cdot\nabla u)\,dx\\
&=-\frac{d}{dt}\int_\R b\cdot\nabla u\cdot b \,dx+\int_\R b_t\cdot \nabla u\cdot b\,dx+\int_\R b\cdot\nabla u\cdot b_t\,dx \\
&\quad-\int_\R b\cdot\nabla u\cdot \nabla u\cdot b\,dx-\int_\R u^j\,b\cdot\nabla \partial_{j} u\cdot b  \,dx\\
&=-\frac{d}{dt}\int_\R b\cdot\nabla u\cdot b \,dx+\int_\R (\Delta b-u\cdot\nabla b+b\cdot\nabla u)\cdot \nabla u\cdot b\,dx\\
&\quad+\int_\R b\cdot\nabla u\cdot (\Delta b-u\cdot\nabla b+b\cdot \nabla u)\,dx-\int_\R b\cdot\nabla u\cdot \nabla u\cdot b\,dx \\
&\quad+\int_\R u\cdot\nabla b\cdot\nabla  u\cdot  b\,dx+\int_\R b\cdot\nabla  u \cdot(u\cdot \nabla b)\,dx\\
&\le -\frac{d}{dt}\int_\R b\cdot\nabla u\cdot b \,dx+\delta\,\|\Delta b\|_{L^2}^2+C(\delta)\left(\|\nabla u\|_{L^3}^3+\|b\|_{L^6}^6\right).
\end{aligned}\end{equation}

Plugging \eqref{3.28}, \eqref{C416}, and \eqref{p2} into \eqref{3.6}  gives rise to
\begin{equation}\begin{aligned}\label{3.133}
&\frac{d}{dt}\int_\R \left(\mu(\rho)|Du|^2+b\cdot\nabla u\cdot b \right)\, dx+\int_\R\rho|\dot u|^2\,dx\\
&\leq\delta\,\|\Delta b\|_{L^2}^2+C(\delta)\left(\|\rho\dot u-b\cdot\nabla b\|_{L^2}\|\nabla u\|_{L^2}^2+\|\nabla u\|_{L^3}^3+\|b\|_{L^6}^6\right).
\end{aligned}\end{equation}

Furthermore, taking the $L^2$ inner product of \eqref{C101}$_3$ with $\Delta b$ and using $\mathrm{div}u=0$ yield
\begin{equation}\begin{aligned}\label{p3}
&\frac{1}{2}\,\frac{d}{dt}\int_\R|\nabla b|^2\, dx+\int_\R|\Delta b|^2\, dx \\
& =\int_\R u\cdot \nabla b\cdot\Delta b\,dx-\int_\R b\cdot\nabla u\cdot\Delta b\,dx \\
&=-\int_\R\partial_i u\cdot \nabla b\cdot\partial_i b\,dx-\int_\R u\cdot\nabla\partial_i b\cdot\partial_i b\,dx-\int_\R b\cdot\nabla u\cdot\Delta b\,dx \\
& \le C\,\|\nabla b\|_{L^2}^{4/3}\|\nabla^2 b\|_{L^2}^{2/3}\|\nabla u\|_{L^3}+\|\Delta b\|_{L^2}\|b\|_{L^6}\|\nabla u\|_{L^3}\\
& \le \delta \,\|\Delta b\|_{L^2}^2+C(\delta)\left(\|\nabla u\|_{L^3}^3+\|\nabla b\|_{L^2}^{4}+\|b\|_{L^6}^6\right).
\end{aligned}\end{equation}

Therefore, combining \eqref{3.133} with \eqref{p3} and   choosing suitably small $\delta$,  there exists a positive constant $\varepsilon_1$ such that if $\|\mu(\rho_0)-1\|_{L^\infty}\le \varepsilon_1$,
\begin{equation}\begin{aligned}\label{3.13}
&\frac{d}{dt}\int_\R \left( 2\,\mu(\rho)|Du|^2+2\,b\cdot\nabla u\cdot b+|\nabla b|^2\right)\, dx+\int_\R\left(2\,\rho|\dot u|^2+|\Delta b|^2\right)\,dx\\
&\leq C\left(\|\rho\dot u-b\cdot\nabla b\|_{L^2}\|\nabla u\|_{L^2}^2+\|\nabla u\|_{L^3}^3+\|b\|_{L^6}^6+\|\nabla b\|_{L^2}^{4}\right)\\
&\leq C\left(\|\rho\dot u-b\cdot\nabla b\|_{L^2}\|\nabla u\|_{L^2}^2+\|b\|_{L^2}^2\|\nabla b\|_{L^2}^4+\|\nabla b\|_{L^2}^{4}\right)\\
&\leq \|\sqrt\rho \dot u\|_{L^2}^2+C\left(\||\nabla b|\,|b|\|_{L^2}^2+\|\nabla u\|_{L^2}^4+\|\nabla b\|_{L^2}^{4}\right),
\end{aligned}\end{equation}
due to Lemma \ref{pl3} and  \eqref{3.7}--\eqref{3.8}. In addition, observe that \eqref{1.5} and \eqref{3.8} imply 
\begin{equation}\begin{aligned}\label{p12}
&\int_\R \left(2\, \mu(\rho)|Du|^2+2\,b\cdot\nabla u\cdot b+|\nabla b|^2\right)\, dx\geq \frac{\underline\mu}{2}\|\nabla u\|_{L^2}^2+\|\nabla b\|_{L^2}^{2}-C\,\|b\|_{L^4}^4,\\
&\int_\R \left(2\, \mu(\rho)|Du|^2+2\,b\cdot\nabla u\cdot b+|\nabla b|^2\right)\, dx\le C\left(\|\nabla u\|_{L^2}^2+\|\nabla b\|_{L^2}^{2}\right).
\end{aligned}\end{equation}
Then multiplying \eqref{3.13} by $t^i$ $(i=0,1)$, it follows from  Gronwall's inequality, \eqref{p12}, \eqref{3.8},  and  \eqref{q166}  that \eqref{3.9} holds.  This finishes the proof of Lemma \ref{le3}.
\end{proof}

\begin{lemma}\label{le4}
Under the conditions of Theorem \ref{T1}, suppose that 
$$\|\mu(\rho_0)-1\|_{L^\infty}\leq \varepsilon_2,$$
with $\varepsilon_2$ ($\le \varepsilon_1$) depending only on  $\bar\rho$ and $\mu$, then there exists a positive constant $C$ depending on $\bar\rho$, $\mu$, $\|\rho_0\|_{L^1}$, $\|\sqrt{\rho_0}u_0\|_{L^2}$, $\|\nabla u_0\|_{L^2}$, and $\|b_0\|_{H^1}$  such that   for $i=1,2$, 
\begin{equation}\label{3.10}
\sup_{t\in[0,T]}t^i\left(\|\sqrt\rho\dot u\|_{L^2}^2+\||\nabla b|\,|b|\|_{L^2}^2\right)+\int_0^T t^i\left(\|\nabla\dot u\|_{L^2}^2+\||\Delta b|\,|b|\|_{L^2}^2\right)\,dt\leq C.
\end{equation}
\end{lemma}
\begin{proof}
We first get, by applying the operator ``$\partial_t+
u\cdot\nabla $" to  \eqref{C101}$_2$,  that
\begin{equation}\begin{aligned}\notag
&\rho\dot u_t+\rho u\cdot\nabla\dot u-\mathrm{div}(2\mu(\rho)D\dot u)+\partial_i\left(\mu(\rho)\partial_i u\cdot\nabla u\right)+\partial_i\left(\mu(\rho)\nabla u^i\cdot\nabla u\right)\\
&+2\,\partial_i\left(\mu(\rho)\nabla u\cdot\nabla u^i\right)+\nabla P_t+u\cdot\nabla\left(\nabla P\right)-\partial_i\left(b^i b\right)_t-\partial_i\left( u^ib\cdot\nabla b \right)=0.
\end{aligned}\end{equation}
Multiplying the above equality by $\dot u$ and integrating the resulting equality by parts, we deduce from \eqref{1.5} and \eqref{p1} that  for any $\delta\in(0,1)$,
\begin{equation}\label{3.34}\begin{aligned}
&\frac{1}{2}\,\frac{d}{dt}\int_\R\rho|\dot{u}|^{2}\,dx+2\int_\R\mu(\rho)\,|D\dot u|^2\,dx\\
&=\int_\R\mu(\rho)\left(\partial_i u\cdot\nabla u\cdot\partial_i\dot u+\nabla u^i\cdot\nabla u\cdot\partial_i\dot u+2\,\partial_i\dot u\cdot\nabla u\cdot\nabla u^i\right)\,dx\\
&\quad+\int_\R \left(P_t+u\cdot\nabla P\right)\mathrm{div}\dot u\,dx+\int_\R \dot u\cdot\nabla u\cdot\nabla P\,dx\\
&\quad-\int_\R\partial_i\dot u\cdot\left(b^i b\right)_t\,dx-\int_\R u\cdot\nabla\dot u\cdot( b\cdot\nabla b)\,dx\\
&\leq \delta\,\|\nabla\dot u\|_{L^2}^2+C(\delta)\|\nabla u\|_{L^4}^4+\int_\R \left(P_t+u\cdot\nabla P\right)\partial_i u\cdot\nabla u^i\,dx\\
&\quad-\int_\R P\,\partial_i u\cdot\nabla  \dot u^i\,dx-\int_\R\partial_i\dot u\cdot\left(b^i b\right)_t\,dx-\int_\R u\cdot\nabla\dot u\cdot( b\cdot\nabla b)\,dx.
\end{aligned}
\end{equation}

By using integration by parts and $\mathrm{div} u=0$, we have 
\begin{equation}\begin{aligned}\notag
&\int_\R \left(P_t+u\cdot\nabla P\right)\partial_i u\cdot\nabla u^i\,dx-\int_\R P\,\partial_i u\cdot\nabla  \dot u^i\,dx\\
&=\frac{d}{dt}\int_\R P \,\partial_i u\cdot\nabla u^i\, dx-2\int_\R P\,\partial_i u\cdot\nabla \left(\dot u^i-u\cdot\nabla u^i\right)\,dx\\
&\quad-\int_\R P\,  u\cdot\nabla\left(\partial_i u\cdot\nabla u^i\right)\,dx-\int_\R P\,\partial_i u\cdot\nabla  \dot u^i\,dx\\
&=\frac{d}{dt}\int_\R P \,\partial_i u\cdot\nabla u^i\, dx-3\int_\R P\, \partial_i u\cdot\nabla \dot u^i\,dx+2\int_\R P\, \partial_i u\cdot\nabla u\cdot\nabla u^i\,dx\\
&\quad +2\int_\R P\, u\cdot\nabla \partial_i u\cdot\nabla u^i\,dx-\int_\R P\, u\cdot\nabla\left(\partial_i u\cdot\nabla u^i\right)\,dx\\
&=\frac{d}{dt}\int_\R P\, \partial_i u\cdot\nabla u^i\,dx-3\int_\R P\, \partial_i u\cdot\nabla \dot u^i\,dx+2\int_\R P\,\partial_i u\cdot\nabla u\cdot\nabla u^i\,dx\\
&\le \frac{d}{dt}\int_\R P\, \partial_i u\cdot\nabla u^i\,dx+\delta\,\|\nabla\dot u\|_{L^2}^2+C(\delta)\left(\|\nabla u\|_{L^4}^4+\|P\|_{L^4}^4\right).
\end{aligned}
\end{equation}
Similarly, by virtue of \eqref{C101}$_3$ and $\mathrm{div} b=0$, we obtain
\begin{equation}\begin{aligned}\notag
&-\int_\R\partial_i\dot u\cdot\left(b^i b\right)_t\,dx\\
&=\int_\R \left(u\cdot\nabla b-b\cdot\nabla u-\Delta b\right)\cdot\nabla\dot u\cdot b\, dx+\int_\R b\cdot\nabla\dot u\cdot \left(u\cdot\nabla b-b\cdot\nabla u-\Delta b\right)\, dx\\
&\le\int_\R \left(u\cdot\nabla b\cdot\nabla\dot u\cdot b+ b\cdot\nabla\dot u\cdot (u\cdot\nabla b)\right)\, dx+2\int_\R\left(|b|^2\,|\nabla u|\,|\nabla\dot u|+|b|\,|\Delta b|\,|\nabla\dot u|\right)\,dx\\
&=-\int_\R u\cdot\nabla \partial_i\dot u\cdot b\, b^i\, dx+2\int_\R\left(|b|^2\,|\nabla u|\,|\nabla\dot u|+|b|\,|\Delta b|\,|\nabla\dot u|\right)\,dx\\
&=\int_\R \left( u\cdot\nabla\dot u\cdot( b\cdot\nabla b)+ b\cdot\nabla u\cdot\nabla\dot u\cdot b\right)\,dx+2\int_\R\left(|b|^2\,|\nabla u|\,|\nabla\dot u|+|b|\,|\Delta b|\,|\nabla\dot u|\right)\,dx\\
&\le\int_\R  u\cdot\nabla\dot u\cdot( b\cdot\nabla b)\, dx+\delta\,\|\nabla\dot u\|_{L^2}^2+C(\delta)\left(\|\nabla u\|_{L^4}^4+\|b\|_{L^8}^8+\||\Delta b|\,|b|\|_{L^2}^2\right).
\end{aligned}
\end{equation}

Inserting the above estimates into \eqref{3.34} and choosing suitable $\delta$ give rise to
\begin{equation}\label{3.37}\begin{aligned}
&\frac{d}{dt}\int_\R\left(\rho|\dot{ u }|^{2}-2P\,\partial_i u\cdot\nabla u^i\right)\,dx+\underline\mu\int_\R|\nabla\dot u|^2\,dx\\
&\leq C\,\left(\|\nabla u\|_{L^4}^4+\|P\|_{L^4}^4+\|b\|_{L^8}^8\right)+C'\,\||\Delta b|\,|b|\|_{L^2}^2,
\end{aligned}
\end{equation}
where we have used \eqref{1.5} and the following fact owing to \eqref{p1}:
\begin{equation}\begin{aligned}\notag
2\int_\R|D\dot u|^2\,dx 
&=\int_\R\left(|\nabla \dot u|^2+(\mathrm{div}\dot u)^2\right)\,dx 
\geq\int_\R|\nabla \dot u|^2\,dx-\int_\R|\nabla u|^4\,dx.
\end{aligned}\end{equation}

Next, we adopt the method in \cite{lvhuang,lsx} to handle the magnetic field. Precisely, for $a_1,a_2\in\{-1,0,1\}$, set 
$$B(a_1,a_2)=a_1b^1+a_2b^2,\quad U(a_1,a_2)=a_1u^1+a_2u^2.$$
Then from \eqref{C101}$_3$, it holds that
\begin{equation}\notag
B_t+u\cdot\nabla B-b\cdot\nabla U-\Delta B=0.
\end{equation}
Multiplying this equality by $4\Delta |B|^2 B$ and integrating the result by parts show that
\begin{equation}\label{p6}\begin{aligned}
&\frac{d}{dt}\int_\R|\nabla|B|^2|^2\,dx+2\int_\R|\Delta |B|^2|^2\,dx\\
&=2\int_\R u\cdot\nabla|B|^2\Delta |B|^2\,dx-4\int_\R b\cdot\nabla U\cdot B \Delta |B|^2\,dx+4\int_\R|\nabla B|^2\Delta |B|^2\,dx\\
&\le C\left(\|\nabla u\|_{L^4}\|\nabla B\|_{L^4}^2\|B\|_{L^8}^2+\|\nabla U\|_{L^4}\|b\|_{L^8}^2\|\Delta |B|^2\|_{L^2}+\|\nabla B\|_{L^4}^2\|\Delta |B|^2\|_{L^2}\right)\\
&\le \|\Delta |B|^2\|_{L^2}^2+C\left(\|\nabla u\|_{L^4}^4+\|\nabla b\|_{L^4}^4+\|b\|_{L^8}^8\right).
\end{aligned}
\end{equation}
Moreover, denote 
$$\Phi(t)\triangleq \|\nabla|B(1,0)|^2\|_{L^2}^2+\|\nabla|B(0,1)|^2\|_{L^2}^2+\|\nabla|B(1,1)|^2\|_{L^2}^2+\|\nabla|B(1,-1)|^2\|_{L^2}^2.$$
Direct calculation infers that
\begin{equation}\label{p9}
\||\nabla b|\,|b|\|_{L^2}^2\le \Phi(t)\le C\,\||\nabla b|\,|b|\|_{L^2}^2,
\end{equation}
\begin{equation}\label{p7}\begin{aligned}
\||\Delta b|\,|b|\|_{L^2}^2&\le C\,\|\nabla b\|_{L^4}^4+\|\Delta|B(1,0)|^2\|_{L^2}^2+\|\Delta|B(0,1)|^2\|_{L^2}^2\\
&\quad+\|\Delta|B(1,1)|^2\|_{L^2}^2+\|\Delta|B(1,-1)|^2\|_{L^2}^2.
\end{aligned}
\end{equation}
Therefore, from \eqref{p6} and \eqref{p7}, we arrive at
\begin{equation}\label{p8}\begin{aligned}
\Phi'(t)+\||\Delta b|\,|b|\|_{L^2}^2\le C\left(\|\nabla u\|_{L^4}^4+\|\nabla b\|_{L^4}^4+\|b\|_{L^8}^8\right).
\end{aligned}
\end{equation}

In view of Lemma \ref{pl3}, \eqref{C225}, \eqref{3.7},  and \eqref{1.5}, we infer that for some constant $\varepsilon_2\,(\le \varepsilon_1)$, if $\|\mu(\rho_0)-1\|_{L^\infty}\le \varepsilon_2$,
\begin{equation}\begin{aligned}\label{4}
&\|\nabla u\|_{L^4}^4+\|P\|_{L^4}^4\\
&\leq C\,\|\nabla u\|_{L^2}^2\|\rho\dot u-b\cdot\nabla b\|_{L^{2}}^2+C\,\|\Delta^{-1}\mathrm{div}(\rho\dot u-\mathrm{div}(b\otimes b))\|_{L^4}^4\\
&\leq C\,\|\nabla u\|_{L^2}^2\left(\|\sqrt\rho\dot u\|_{L^2}^2+\||\nabla b|\,|b|\|_{L^2}^2\right)+C\,\|\rho\dot u\|_{L^{4/3}}^4+C\,\|b\|_{L^8}^8\\
&\leq  C\,\|\nabla u\|_{L^2}^2\left(\|\sqrt\rho\dot u\|_{L^2}^2+\||\nabla b|\,|b|\|_{L^2}^2\right)+C\,\|\sqrt\rho\dot u\|_{L^{2}}^4+C\,\|b\|_{L^8}^8.
\end{aligned}\end{equation}
By \eqref{3.8}, we  have
\begin{equation}\label{p4}
\|b\|_{L^8}^8=\||b|^2\|_{L^4}^4\le C\,\|b\|_{L^4}^4\||\nabla b|\,|b|\|_{L^2}^2\le C\,\|\nabla b\|_{L^2}^2\||\nabla b|\,|b|\|_{L^2}^2.
\end{equation} 
In addition, note that \eqref{p1},  \eqref{C416},  \eqref{3.7}, and Lemma \ref{pl3} imply that if $\|\mu(\rho_0)-1\|_{L^\infty}\le \varepsilon_2$,
\begin{equation}\label{3.38}\begin{aligned}
\left|\int_\R P\,\partial_i u\cdot\nabla u^i\,dx\right|
&\leq C\,\|\rho\dot u-b\cdot\nabla b\|_{L^{2}}\|\nabla u\|_{L^2}^2+C\,\|\nabla u\|_{L^3}^3\\
&\leq \frac{1}{4}\left(\|\sqrt\rho\dot u\|_{L^{2}}^2+\||\nabla b|\,|b|\|_{L^2}^2\right)+C\,\|\nabla u\|_{L^2}^4.
\end{aligned}
\end{equation}
Now, adding \eqref{3.37} to \eqref{p8} multiplied by $C'+1$ and applying \eqref{4}, \eqref{p4} conclude
\begin{equation}\begin{aligned}\label{p27}
&\frac{d}{dt}\int_\R\left(\rho|\dot{ u }|^{2}-2P\,\partial_i u\cdot\nabla u^i\right)\,dx+(C'+1)\Phi'(t)+\underline\mu\,\|\nabla\dot u\|_{L^2}^2+\||\Delta b|\,|b|\|_{L^2}^2\\
&\leq C\left(\|\nabla u\|_{L^4}^4+\|\nabla b\|_{L^4}^4+\|P\|_{L^4}^4+\|b\|_{L^8}^8\right)\\
&\leq  C\left(\|\nabla u\|_{L^2}^2+\|\nabla b\|_{L^2}^2+\|\sqrt\rho\dot u\|_{L^{2}}^2\right)\left(\|\sqrt\rho\dot u\|_{L^2}^2+\||\nabla b|\,|b|\|_{L^2}^2\right)\\
&\quad+C\,\|\nabla b\|_{L^2}^2\|\Delta b\|_{L^2}^2.
\end{aligned}
\end{equation}
Multiplying this inequality by $t^i$ $(i=1,2)$ and using Gronwall's inequality, we thus deduce from \eqref{p9}, \eqref{3.38}, \eqref{3.8}, \eqref{q166}, and \eqref{3.9}  
that \eqref{3.10} holds provided $\|\mu(\rho_0)-1\|_{L^\infty}\le \varepsilon_2$. The proof of Lemma \ref{le4} is completed.
\end{proof}

Next,  we provide the additional integrability of $u$ and $b$ in the following lemmas, which play  important roles in the later analysis.
\begin{lemma}\label{le10}
Under the conditions of Theorem \ref{T1}, for any $p\in(1,2)$, there exists a positive constant $C$ depending on $p$, $\bar\rho$, $\mu$,  $\|\sqrt{\rho_0}u_0\|_{L^2}$, and $\|b_0\|_{L^1\cap L^2}$ such that 
\begin{equation}\label{p22}
\|b\|_{L^p} \le C,\quad \|b\|_{L^2}^{2}\le C\,t^{-(2-p)/p},\quad t\in(0,T].
\end{equation}
\end{lemma}
\begin{proof}
For any $p\in(1,2)$, taking the inner product of  \eqref{C101}$_3$ with $p|b|^{p-2} b$ yields 
\begin{equation}\begin{aligned}\notag
\frac{d}{dt}\int_\R |b|^p\,dx+p\int_\R |b|^{p-2}|\nabla b|^2\,dx+\frac{4(p-2)}{p}\int_\R |\nabla|b|^{p/2}|^2\,dx
\le p\int_\R |b|^p|\nabla u|\,dx.
\end{aligned}\end{equation}
Notice that $|\nabla|b|^{p/2}|^2\le p^2|b|^{p-2}|\nabla b|^2/4$, then for any $\delta\in(0,1)$,
\begin{equation}\begin{aligned}\notag
\frac{d}{dt}\int_\R|b|^p\,dx+\frac{4(p-1)}{p}\int_\R|\nabla|b|^{p/2}|^2\,dx
&\le p\int_\R|b|^p|\nabla u|\,dx\\
&\le C\,\|\nabla u\|_{L^2}\||b|^{p/2}\|_{L^2}\|\nabla |b|^{p/2}\|_{L^2}\\
&\le \delta\,\|\nabla |b|^{p/2}\|_{L^2}^2+C(\delta)\,\|\nabla u\|_{L^2}^2\|b\|_{L^p}^p.
\end{aligned}\end{equation}
By taking suitably small $\delta$ and applying Gronwall's inequality and \eqref{3.8}, we obtain 
\begin{equation}\label{p23}
\sup_{ t\in[0,T]}\|b\|_{L^p}^p+\int_0^T\|\nabla |b|^{p/2}\|_{L^2}^2\,dt\le C,
\end{equation}
which along with direct calculation shows that
\begin{equation}\label{p21}
\int_0^T\|b\|_{L^2}^{2p/(2-p)}\,dt
\le C\int_0^T \||b|^{p/2}\|_{L^2}^{2p/(2-p)}\|\nabla |b|^{p/2}\|_{L^2}^2\,dt\le C.
\end{equation}

In the similar way, we get 
$$\frac{d}{dt}\|b\|_{L^2}^2+\|\nabla b\|_{L^2}^2\le C\,\|\nabla u\|_{L^2}^2\|b\|_{L^2}^2.$$
Thereby, multiplying it by $t\|b\|_{L^2}^{2p/(2-p)-2}$ yields 
$$\frac{d}{dt}\left(t\,\|b\|_{L^2}^{2p/(2-p)}\right)\le C\,t\,\|\nabla u\|_{L^2}^2\|b\|_{L^2}^{2p/(2-p)}+\|b\|_{L^2}^{2p/(2-p)}.$$
Utilizing Gronwall's inequality again together with  \eqref{3.8} and \eqref{p21} leads to that $t\,\|b\|_{L^2}^{2p/(2-p)}\le C$, which and \eqref{p23} give \eqref{p22}.
 The proof of Lemma \ref{le10} is finished.
\end{proof}

To obtain the integrability of $u$, we first adopt the same methods as in \cite[Lemmas 3.3-3.4]{FLW} to establish  the following two  cancellation mechanisms.
\begin{lemma}\label{le0}
The following conclusions hold
\begin{equation}\label{cancel}
\int_\R\rho\dot u(x,t)\,dx=0,~~\mathrm{for}~\mathrm{a.e.}~t\in(0,T],
\end{equation}
\begin{equation}\label{cancel1}
\int_\R\rho u(x,t)\,dx=\int_\R\rho_0 u_0(x)dx=0,~~\mathrm{for}~\mathrm{a.e.}~t\in(0,T].
\end{equation}
\end{lemma}
\begin{proof}
Set $\psi_r$ be the  standard truncation function satisfying $0\le \psi_r\in C_0^\infty(B_{2r})$, $\psi_r\equiv  1$ in $B_{r}$, $|\nabla^i\psi_r|\le Cr^{−i}$, $i = 1, 2$. For any fixed constant vector $a=(a_1,a_2)$ and $\nabla^\bot=(-\partial_2,\partial_1)$,  we define
\[
f(x)=a_2x_1-a_1x_2,\quad 
\Psi_r=\nabla^\perp(f\psi_r)
=\psi_r a+f \nabla^\perp\psi_r.
\]
Then  $\|\Psi_r\|_{L^\infty}\leq C|a|$, $\mathrm{supp}\nabla\Psi_{r}\subset B_{2r}\setminus B_{r}$, $\|\nabla\Psi_r\|_{L^\infty}\leq C|a|/r$, and $\|\nabla\Psi_r\|_{L^2}\leq C|a|$.

For any $\phi(t)\in C_0^\infty([0,T])$, multiplying  \eqref{C101}$_2$ by $\Psi_r\,\phi$ and integrating the resulting equality by parts, we obtain after using  $\mathrm{div}\Psi_r=0$ and \eqref{1.7} that
\begin{equation}\begin{aligned}\notag
&\left|\int_0^T\int_\R\rho\dot u\cdot\Psi_r\,\phi\,dxdt\right|\\
&=\left|-2\int_0^T\int_{B_{2r}\setminus B_{r}}\mu(\rho)Du:\nabla\Psi_r\,\phi\,dxdt-\int_0^T\int_{B_{2r}\setminus B_{r}}b\otimes b:\nabla\Psi_r\,\phi\,dxdt\right|\\
&\le C|a|\int_0^T\|\nabla u\|_{L^2(B_{2r}\setminus B_{r})}\,dt+\frac{C|a|}{r}\int_0^T\|b\|_{L^2}^2\,dt\quad\rightarrow0
~~\mathrm{as}~~r\rightarrow\infty.
\end{aligned}\end{equation}
Observe that $\|\rho\dot u\|_{L^1}\le M^{1/2}\|\sqrt\rho\dot u\|_{L^2}\in L^1(0,T)$ and $\Psi_r\to a$ pointwise as $r\rightarrow\infty$, thereby using the dominated convergence theorem in above equality leads to
\[
a\cdot\int_0^T\left(\int_\R\rho\dot u\,dx\right)\phi(t)\,dt=0.
\]
Due to the arbitrariness of $a\in\R$ and $\phi\in C_0^\infty([0,T])$, we derive \eqref{cancel} immediately.

Now, we take $\phi(t)\in C_0^\infty([0,T))$ and test \eqref{C101}$_2$ against $\Psi_r\,\phi$ again to get
\begin{equation}\label{CC403}
\begin{split}
&\int_0^T\phi'(t)\int_\R\rho u\cdot\Psi_r\,dx\,dt+\phi(0)\int_\R\rho_0u_0\cdot\Psi_r\,dx\\
&=-\int_0^T\int_\R
\bigl(\rho u\otimes u-2\mu(\rho)Du-b\otimes b\bigr):\nabla\Psi_r\,\phi\,dx\,dt.
\end{split}
\end{equation}
where the right hand side  satisfies
\begin{equation*}
\begin{split}
&\left|\int_0^T\int_\R
\bigl(\rho u\otimes u-2\mu(\rho)Du-b\otimes b\bigr):\nabla\Psi_r\,\phi\,dx\,dt\right|\\
&\leq\frac{C|a|}{r}
\int_0^T\left(\|\sqrt\rho u\|_{L^2}^2+\|b\|_{L^2}^2\right)\,dt+ C|a|
\int_0^T\|\nabla u\|_{L^2(B_{2r}\setminus B_r)}\,dt\\
& \rightarrow 0 ~~\mathrm{as}~~r\rightarrow\infty.
\end{split}
\end{equation*}
Applying the dominated convergence theorem in  \eqref{CC403}, where $\|\rho u\|_{L^1}\le M^{1/2}\|\sqrt\rho u\|_{L^2}$, then yields
\begin{equation}\notag
a\cdot\int_0^T\left(\int_\R\rho u\,dx\right)\phi'(t)\,dt+\phi(0)a\cdot\,
\left(\int_\R\rho_0u_0\,dx\right)=0,
\end{equation}
which combined with the arbitrariness of $a\in\R$ and $\phi\in C_0^\infty([0,T))$ infers \eqref{cancel1}. The proof of Lemma \ref{le0} is finished.
\end{proof}

\begin{lemma}\label{le11}
Under the conditions of Theorem \ref{T1},  suppose that
$$\|\mu(\rho_0)-1\|_{L^\infty}\leq \varepsilon_3,$$
with $\varepsilon_3$ ($\le\varepsilon_2$) depending only on $\bar\rho$ and $\mu$, then for a.e. $t\in(0,T]$, the following states hold.

1. There exists a unique far field velocity $u_\infty(t)$ such that $u(x,t)\rightarrow u_\infty(t)$ as $|x|\rightarrow\infty$.

2. For any  $0<\alpha<\beta<1/2$, there exists a positive  constant $C$  depending on $\alpha$, $\beta$, $\bar\rho$, $\mu$,  $\|\rho_0\|_{L^1}$, $\|\sqrt{\rho_0}u_0\|_{L^2}$, and $\|b_0\|_{L^1\cap L^2}$  such that 
\begin{equation}\label{p10} 
\| u-u_\infty\|_{L^{2/\alpha}}+|u_\infty(t)|\le  C\left(1+\|\rho\,\tilde x^2\|_{L^1}\right)^{\beta/2}\|\sqrt{\rho}\dot u\|_{L^2}+C\,\|\nabla b\|_{L^2}.
\end{equation}
\end{lemma}
\begin{proof}
To begin with, we rewrite \eqref{C101}$_2$ as
\begin{equation}\label{ude}
-\Delta u+\nabla P=\mathrm{div}(2(\mu(\rho)-1)Du)-(\rho\dot u-b\cdot\nabla b).\end{equation}
Then we consider the following two elliptic systems:
\begin{equation}\label{p11}
-\Delta w_1+\nabla P_1=\mathrm{div}(2(\mu(\rho)-1)Dw_1)-\rho\dot u,\quad  \mathrm{div}w_1=0,
\end{equation}
\begin{equation}\label{p13}
-\Delta w_2+\nabla P_2=\mathrm{div}(2(\mu(\rho)-1)Dw_2)+b\cdot\nabla b,\quad  \mathrm{div}w_2=0.
\end{equation}

On the one hand, for any $0<\alpha<\beta<1/2$,  we  define $p=2/(1+\alpha)$, $ p^*=2/\alpha$, and 
$$X_1\triangleq\{w\in L^{p^*}|\mathrm{div}w=0, \nabla w\in  L^2,  \nabla w|x|^\alpha\in L^2\},
\quad
\|w\|_{X_1}=\|\nabla w\|_{L^2}+\|\nabla w\,|x|^\alpha\|_{L^2}.
$$
Then $X_1$ is a Banach space due to Lemma \ref{pl2}. Moreover, we set
$$\mathcal T_1(w_1)\triangleq-\Delta^{-1}\mathbb  P\operatorname{div}(2(\mu-1)Dw_1)+\Delta^{-1}\mathbb{P}(\rho\dot u),
$$
where  the operator is understood with its \(L^{p^*}\) representative and ``$\mathbb{P}\triangleq \mathrm{Id}-\nabla\Delta^{-1}\mathrm{div}$" is the Leray projector.

Since for a.e. $t\in[0,T]$, $\rho\dot u$ satisfies 
$$ \rho\dot u\in L^1\cap L^2,\quad \int_\R \rho\dot u\,dx=0,\quad\int_\R|\rho\dot u|\,|x|\,dx\le \|\rho\tilde x^2\|_{L^1}^{1/2}\|\sqrt\rho\dot u\|_{L^2}<\infty,$$
due to \eqref{1.7} and \eqref{cancel}, it must belong to $\mathcal{H}^1$ according to \cite[Remark2.2]{FLW}. 
Based on this, we adopt the key idea in \cite[Lemma 3.2]{FLW} to estimate $\mathcal T_1(w_1)$. More precisely, for any $0<\alpha<\beta<1/2$, it holds that
$$\|\Delta^{-1}\mathbb  P\operatorname{div}(2(\mu-1)Dw_1)\|_{X_1}\le C_1\,\|\mu(\rho_0)-1\|_{L^\infty}\|w_1\|_{X_1},$$
\begin{equation}\begin{aligned}\label{har}
&\|\Delta^{-1}\mathbb{P}(\rho\dot u)\|_{X_1}+\|\rho\dot u\|_{\mathcal{H}^1}\\
&\le C\left(\|\rho\|_{L^1}^{1/2}+\|\rho\|_{L^\infty}^{1/2}+\|\rho\|_{L^1}^{(1-\beta)/2}\|\rho\,|x|^2\|_{L^1}^{\beta/2}\right)\,\|\sqrt{\rho}\dot u\|_{L^2}\\
&\le C\left(1+\|\rho\,\tilde x^2\|_{L^1}\right)^{\beta/2}\|\sqrt{\rho}\dot u\|_{L^2},
\end{aligned}\end{equation}
owing to \eqref{3.7}. Then,
\begin{equation}\label{C339}
\begin{split}
\|\mathcal T_1(w_1)\|_{X_1}
&\le \frac{1}{2}\|w_1\|_{X_1}+ C\left(1+\|\rho\,\tilde x^2\|_{L^1}\right)^{\beta/2}\|\sqrt{\rho}\dot u\|_{L^2},
\end{split}
\end{equation}
$$
\|\mathcal T_1(w_1)-\mathcal T_1(\widetilde w_1)\|_{X_1}
\le C_1\|\mu(\rho_0)-1\|_{L^\infty}\|w_1-\widetilde w_1\|_{X_1}\le \frac{1}{2}\|w_1-\widetilde w_1\|_{X_1},
$$
provided $\|\mu(\rho_{0})-1\|_{L^\infty}\leq \varepsilon_{31}\triangleq 1/(2C_1)$. Since the \(A_2\) characteristics of the weights \(|x|^{2\alpha}\) are uniformly bounded for \(0<\alpha<1/2\), $\varepsilon_{31}$ can be chosen independently of $\alpha$. Together with the weighted Sobolev inequality in Lemma \ref{pl2} and the choice of the \(L^{p^*}\) representative, \eqref{C339}  shows that \(\mathcal T_1\) maps \(X_1\) into itself.
Hence, with the help of the Banach's fixed point theorem,  $\mathcal{T}_1$ has a unique fixed point, that is, the system \eqref{p11}  admits a unique weak solution  $w_1$ in $X_1$. In particular,  \eqref{C339} combined with Lemma \ref{pl2} gives that
\begin{equation}\begin{aligned}\label{w1}
\|w_1\|_{L^{p^\ast}}
\le C\,\|\nabla w_1\,|x|^\alpha\|_{L^2}\le C\left(1+\|\rho\,\tilde x^2\|_{L^1}\right)^{\beta/2}\|\sqrt{\rho}\dot u\|_{L^2}.
\end{aligned}\end{equation}


On the other hand, note that $b\cdot\nabla b=\mathrm{div}(b\otimes b)$, where $b\otimes b\in L^{2/(1+\alpha)}\cap  L^{2}$ with $\alpha$ as above.  Similarly, we  define the  Banach space $X_2$ and operator $\mathcal{T}_2$,
$$X_2\triangleq\{w\in L^{p^*}|\mathrm{div}w=0,\ \nabla w\in L^p\cap L^2\},
\quad
\|w\|_{X_2}=\|\nabla w\|_{L^p}+\|\nabla w\|_{L^2},
$$
$$\mathcal T_2(w_2)\triangleq-\Delta^{-1}\mathbb  P\operatorname{div}(2(\mu-1)Dw_2+b\otimes b).
$$
The operator is understood with its \(L^{p^*}\) representative. By the homogeneous Sobolev inequality, \(\mathcal T_2\) maps \(X_2\) into itself. 
By the boundedness of ``$\nabla\Delta^{-1}\mathbb{P}\mathrm{div}$" in $L^r$ $(1<r<\infty)$ and \eqref{3.7}, we get 
\begin{equation}\begin{aligned}\label{p15}
\|\mathcal T_2(w_2)\|_{X_2}
&\le C_2\,\|\mu(\rho_0)-1\|_{L^\infty}\| w_2\|_{X_{2}}+C\left(\|b\|_{L^{4/(1+\alpha)}}^2+\|b\|_{L^4}^2\right)\\
&\le \frac{1}{2}\|w_2\|_{X_2}+C\left(\|b\|_{L^{4/(1+\alpha)}}^2+\|b\|_{L^4}^2\right),
\end{aligned}\end{equation}
$$
\|\mathcal T_2(w_2)-\mathcal T_2(\widetilde w_2)\|_{X_2}
\le C_2\,\|\mu(\rho_0)-1\|_{L^\infty}\|w_2-\widetilde w_2\|_{X_2}\le \frac{1}{2}\|w_2-\widetilde w_2\|_{X_2},
$$
provided $\|\mu(\rho_0)-1\|_{L^\infty}\le \varepsilon_{32}\triangleq 1/(2C_2)$, where $\varepsilon_{32}$ is independent of $\alpha$ due to $p\in(4/3,2)$. Applying Banach's fixed point theorem,  the system \eqref{p13} admits a unique weak solution  $w_2$ in $X_2$. By virtue of  \eqref{p15}, \eqref{3.8}, and \eqref{p22}, we have 
\begin{equation}\label{p16}
\|w_2\|_{L^{p^\ast}}\le C\,\|\nabla w_2\|_{L^{p}}\le C\left(\|b\|_{L^{2/(1+\alpha)}}\|\nabla b\|_{L^2}+\|b\|_{L^2}\|\nabla b\|_{L^2}\right)\le  C\,\|\nabla b\|_{L^2}.
\end{equation}
 
Furthermore, we set $v\triangleq w_1+w_2$. \eqref{w1}  as well as \eqref{p16} infers that
\begin{equation}\label{vd}
\| v\|_{L^{2/\alpha}}\le  C\left(1+\|\rho\,\tilde x^2\|_{L^1}\right)^{\beta/2}\|\sqrt{\rho}\dot u\|_{L^2}+C\,\|\nabla b\|_{L^2}.
\end{equation}
We deduce from \eqref{ude}--\eqref{p13} that for $z=u-v$, $Q=P-P_1-P_2$, 
$$-\mathrm{div}(2\mu(\rho)Dz)+\nabla Q=0,\quad \mathrm{div} z=0,\quad \nabla z\in L^2.$$
In view of the density of  $C_{0,0}^\infty$ functions in the space $D^{1}_{0,0}$ containing all equivalence classes up to a constant, we can test $z$ in the above equation to derive $\int_\R2\mu|Dz|^2\,dx=0$, which together with \eqref{1.5} and $2\|Dz\|_{L^2}^2=\|\nabla z\|_{L^2}^2$  implies $\nabla z=\nabla (u-v)=0$.
 Therefore, there is a constant vector field $u_\infty(t)$ such that  $u=v+u_\infty$.

In addition,  it follows from \eqref{cancel1} and \eqref{3.7} that
\begin{equation}\notag
\int_\R\rho v \,dx=\int_\R\rho u\,dx-u_\infty\int_\R\rho \,dx=-u_\infty M,
\end{equation}
which gives $u_\infty(t)=-\frac{1}{M}\int_\R\rho v(x,t) \,dx$ and 
\begin{equation}\label{ui}
|u_\infty(t)|\le \frac{1}{M}\|\rho\|_{L^{2/(2-\alpha)}}\|v\|_{2/\alpha}\le C\left(1+\|\rho\,\tilde x^2\|_{L^1}\right)^{\beta/2}\|\sqrt{\rho}\dot u\|_{L^2}+C\,\|\nabla b\|_{L^2},
\end{equation}
owing to \eqref{3.7} and \eqref{vd}.

Finally, we deduce from \eqref{1.7} that for a.e. $t\in[0,T]$, $\nabla v=\nabla u\in H^1$. Then $v\in W^{1,2/\alpha}$. The standard Sobolev embedding theorem indicates that $v\in C^{1-\alpha}$ and $v\rightarrow 0$ as $|x|\rightarrow \infty$, which is guaranteed by the  uniform continuousness of $v$ in $\mathbb{R}^2$ and $v\in L^{2/\alpha}$. We thus conclude Lemma \ref{le11}  from this fact,  \eqref{vd}, and \eqref{ui} provided $\|\mu(\rho_0)-1\|_{L^\infty}\le \varepsilon_3\triangleq\min\{\varepsilon_{31},\varepsilon_{32},\varepsilon_2\}$.
\end{proof}

Next, we recall the following non-uniform estimates on the  density, whose proof are showed in \cite{FLW,lx}.
\begin{lemma}\label{le9}
Under the conditions of Theorem \ref{T1}, there exists a positive constant $C$ depending on $\|\rho_0\,\tilde x^2\|_{L^1}$,  $\|\sqrt{\rho_0}u_0\|_{L^2}$, and $\|b_0\|_{L^2}$ such that 
\begin{equation}\label{k37} 
\int_\R \rho(x,t)\,\tilde x^2\,dx\le C(1+t)^2,\quad t\in(0,T].
\end{equation}
In particular, there exists a positive constant $N$ depending  on $\|\rho_0\|_{L^1}$, $\|\rho_0\,\tilde x^2\|_{L^1}$, $\|\sqrt{\rho_0}u_0\|_{L^2}$, and $\|b_0\|_{L^2}$  such that 
\begin{equation}\label{k38} 
\int_{B_{N(1+t)}}\rho(x,t)\,dx\geq \frac{1}{2}\,\|\rho_0\|_{L^1},\quad t\in(0,T].
\end{equation}
\end{lemma}

As a consequence of Lemmas \ref{le10}--\ref{le9}, we can obtain faster decay rates of the solutions than Lemmas \ref{le2}--\ref{le4}. 
\begin{lemma}\label{le5} 
Under the conditions of Theorem \ref{T1},   suppose that
$$\|\mu(\rho_0)-1\|_{L^\infty}\leq \varepsilon_3,$$
with $\varepsilon_3$ given by Lemma \ref{le11}, then 
for any  $\eta\in(0,1)$, there exists a positive  constant $C$ depending on $\eta$, $\bar\rho$, $\mu$, $\|\rho_0\|_{L^1}$, $\|\rho_0\,\tilde x^2\|_{L^1}$, $\|\sqrt{\rho_0}u_0\|_{L^2}$, $\|\nabla u_0\|_{L^2}$, and $\|b_0\|_{L^1\cap H^1}$ such that
\begin{equation}\begin{aligned}\label{3.11}
&\sup_{t\in[0,T]}t^\eta\left(\|\sqrt\rho u\|_{L^2}^2+\|b\|_{L^2}^2\right)+\int_0^Tt^\eta\left(\| \nabla u\|_{L^{2} }^2+\|\nabla b\|_{L^2}^2\right)\,dt\le C,\\
&\sup_{t\in[0,T]}t^{1+\eta}\left(\|\nabla u\|_{L^2}^2+\|\nabla b\|_{L^2}^2\right)+\int_0^T t^{1+\eta}\left(\|\sqrt\rho\dot u\|_{L^2}^2+\|\Delta b\|_{L^2}^2\right)\,dt\le C,\\
&\sup_{t\in[0,T]}t^{2+\eta}\left(\|\sqrt\rho\dot u\|_{L^2}^2+\||\nabla b||b|\|_{L^2}^2\right)+\! \int_0^Tt^{2+\eta}\left(\|\nabla \dot u\|_{L^2}^2+\||\Delta b||b|\|_{L^2}^2\right)\,dt\leq C.
\end{aligned}\end{equation} 
\end{lemma}
\begin{proof}
First, for any $\eta\in (0,1)$, multiplying \eqref{3.12} by $t^{\eta}$ and integrating the resulting equality over $[0,T]$, we deduce from  \eqref{1.5}, \eqref{3.8}, \eqref{3.9}, \eqref{3.10}, and \eqref{p22} with $p=4/(3+\eta)$  that for $\sigma(T)\triangleq\min\{1,T\}$,
\begin{equation}\begin{aligned}\label{3.122}
&\sup_{t\in[0,T]}t^\eta\left(\|\sqrt\rho u\|_{L^2}^2+\|b\|_{L^2}^2\right)+2\int_0^T t^\eta\left(\underline\mu\,\|\nabla u\|_{L^2}^2+\|\nabla b\|_{L^2}^2\right)\,dt\\
&\le C\int_0^{\sigma(T)} t^{\eta-1}\left(\|\sqrt\rho u\|_{L^2}^2+\|b\|_{L^2}^2\right)\, dt+C\int_{\sigma(T)}^T t^{\eta-1}\left(\|\sqrt\rho u\|_{L^2}^2+\|b\|_{L^2}^2\right)\, dt\\
&\le C\int_0^{\sigma(T)} t^{\eta-1}\,dt+C\int_{\sigma(T)}^T t^{\eta-1}\left(t^{2/3}\|\sqrt\rho\dot u\|_{L^{2}}^2+\|\nabla b\|_{L^2}^2+t^{-(2-p)/p}\right)\,dt\\
&\le C+C\int_{\sigma(T)}^T \left(t^{\eta-2}+t^{\eta-1-(1+\eta)/2}\right)\,dt\le C,
\end{aligned}\end{equation}
where we have used the following fact:
\begin{equation}\label{nu2}\notag
\begin{split}
\int_\R\rho |u|^2\,dx
&=\int_\R\rho |u-u_{\infty}|^2\,dx-\int_\R\rho |u_{\infty}|^2\,dx
\\
&\leq\int_\R\rho|u-u_{\infty}|^2\,dx\\
&\leq \|\rho\|_{L^{4/3}}\|u-u_\infty\|_{L^8}^2\\
&\leq C\left((1+t)^{2/3} \|\sqrt\rho\dot u\|_{L^2}^2+\|\nabla b\|_{L^2}^2\right),
\end{split}
\end{equation}
due to \eqref{cancel1}, \eqref{3.7}, \eqref{p10} with $\alpha=1/4$, $\beta=1/3$, and \eqref{k37}.
 
Then, \eqref{3.122} together with \eqref{p24}, \eqref{q166}, and \eqref{3.8} leads to
\begin{equation}\begin{aligned}\label{p25}
&\sup_{t\in[0,T]}t^{1+\eta}\,\|b\|_{L^4}^4+4\int_0^Tt^{1+\eta}\,\||\nabla b|\,|b|\|_{L^2}^2\,dt\\
&\le C\int_0^Tt^{1+\eta}\,\|\nabla u\|_{L^2}^2\|b\|_{L^4}^4\,dt+C\int_0^Tt^{\eta}\,\|b\|_{L^2}^2\|\nabla b\|_{L^2}^2\, dt\\
&\le C\int_0^Tt^{\eta}\,\|\nabla u\|_{L^2}^2\,dt+C\int_0^T t^{\eta}\,\|\nabla b\|_{L^2}^2\, dt\le C.
\end{aligned}\end{equation}
Next, multiplying \eqref{3.13} by $t^{1+\eta}$ and  using \eqref{p12}, \eqref{3.9}, \eqref{3.122}, and \eqref{p25} yield
\begin{equation}\begin{aligned}\label{p26}
&\sup_{t\in[0,T]}t^{1+\eta}\left(\|\nabla u\|_{L^2}^2+\|\nabla b\|_{L^2}^2\right)+\int_0^Tt^{1+\eta}\left(\|\sqrt\rho \dot u\|_{L^2}^2+\|\Delta b\|_{L^2}^2\right)\,dt\\
&\leq C\sup_{t\in[0,T]}t^{1+\eta}\,\|b\|_{L^4}^4+C\int_0^Tt^{1+\eta}\left(\||\nabla b|\,|b|\|_{L^2}^2+\|\nabla u\|_{L^2}^4+\|\nabla b\|_{L^2}^{4}\right)\,dt\\
&\quad+C\int_0^Tt^{\eta}\left(\|\nabla u\|_{L^2}^2+\|\nabla b\|_{L^2}^{2}\right)\,dt\le C.
\end{aligned}\end{equation}
Similarly, by virtue of \eqref{p27}, \eqref{3.38}, \eqref{p9}, \eqref{3.9}, \eqref{3.10}, and \eqref{3.122}--\eqref{p26}, we have
\begin{equation}\begin{aligned}\notag
\sup_{t\in[0,T]}t^{2+\eta}\left(\|\sqrt\rho \dot{ u }\|_{L^2}^{2}+\||\nabla b|\,|b|\|_{L^2}^2\right)+\int_0^Tt^{2+\eta}\left(\|\nabla\dot u\|_{L^2}^2+\||\Delta b|\,|b|\|_{L^2}^2\right)\,dt\le C.
\end{aligned}\end{equation}
This combined with \eqref{3.122} and \eqref{p26} finishes the proof of Lemma \ref{le5}.
\end{proof}

\begin{lemma}\label{le6}
Under the conditions of Theorem \ref{T1}, there exist  positive constants $\varepsilon_4$ ($\le \varepsilon_3$)  and $C$ depending on  $\bar\rho$, $\mu$, $q$, 
and the initial values such that if
$$\|\mu(\rho_0)-1\|_{L^\infty}\leq\varepsilon_4,$$
 it holds that
\begin{equation}\label{3.19}
\sup_{t\in[0,T]}\|\rho\|_{H^1\cap W^{1,q}}+\int_0^T\|\nabla u\|_{L^\infty}\,dt\leq C.
\end{equation}
\end{lemma}
\begin{proof}
 First, differentiating  \eqref{C101}$_1$ yields
\begin{equation}\notag
\left(\nabla \rho\right)_t + u \cdot \nabla (\nabla\rho) + \nabla u \cdot \nabla  \rho=0,
\end{equation}
which multiplied by $|\nabla\rho|^{q-2}\nabla\rho$ infers
\begin{equation}\label{3.39}
\frac{d}{dt}\|\nabla \rho \|_{L^q}\leq \|\nabla \rho\|_{L^q}\|\nabla u\|_{L^\infty}.
\end{equation}

Now, we need to control $\|\nabla u\|_{L^\infty}$. To do this, we apply the Leray projector  ``$\mathbb{P}=\mathrm{Id}-\nabla\Delta^{-1}\mathrm{div}$" to \eqref{ude} and derive 
\begin{equation}\label{C422}
\begin{split}
\nabla u
&=-\nabla\Delta^{-1}\mathbb{P}\,\mathrm{div}(2(\mu(\rho)-1)Du)+\nabla\Delta^{-1}\mathbb{P}(\rho\dot u-b\cdot\nabla b)\\
&=\int_\R {K}(x-y)\, \mathrm{div}(2(\mu(\rho)-1)Du)(y)\,dy+\nabla\Delta^{-1}\mathbb{P}(\rho\dot u-b\cdot\nabla b)\\
&\triangleq \mathbf{I}_1(x)+\mathbf{I}_2(x),
\end{split}
\end{equation}
where the kernel $K$ satisfies 
\begin{equation}\label{k2}
|\nabla^m K(x)|\le C(m)|x|^{-(m+1)},\quad \mathrm{for}~\mathrm{any}~ x\in\R\setminus \{0\},\, m\in \mathbb{N}.
\end{equation}

To bound $\mathbf{I}_1(x)$, let us introduce the truncation function $\varphi\in C_0^\infty(\R)$ defined as
\begin{equation}\label{var}
\varphi(x)=
\begin{cases}
1,\quad |x|\leq 1,\\
0,\quad |x|\geq 2,\end{cases}
 0\leq\varphi(x)\leq 1,\quad |\nabla\varphi(x)|\leq 2,\quad\text{for any}\, x\in\R,
\end{equation}
and set $\varphi_\lambda(\cdot)=\varphi\big(\cdot/\lambda\big)$ for some $\lambda=\lambda(t)$ determined later, then we split $\mathbf{I}_1(x)$ by
\begin{equation}
\begin{split}\label{Q1}
\mathbf{I}_1(x)
&=\int_\R\varphi_\lambda(x-y) K(x-y)\,\mathrm{div}(2(\mu(\rho)-1)Du)(y)\,dy\\
&\quad+\int_\R\left(1-\varphi_\lambda(x-y)\right) K(x-y)\,\mathrm{div}(2(\mu(\rho)-1)Du)(y)\,dy\\
&\triangleq \mathbf{I}_{11}(x)+\mathbf{I}_{12}(x).
\end{split}
\end{equation}
In view of \eqref{k2}, \eqref{var}, \eqref{3.7}, and \eqref{1.5}, one gets that for  $p\triangleq\min\{(2+q)/2,9/4\}\in(2,q)$,
\begin{equation}\label{C420}
\begin{split}
|\mathbf{I}_{11}(x)|
&\leq C\int_{\{|x-y|\leq2\lambda\}}|x-y|^{-1}|\nabla\mu(\rho)\cdot Du+(\mu(\rho)-1)\cdot\nabla^2 u|(y)\,dy\\
&\leq C\,\lambda^\frac{p-2}{p}\left(\|\nabla\mu(\rho)\cdot Du\|_{L^p}+\|\mu(\rho)-1\|_{L^\infty}\|\nabla^2u\|_{L^p}\right)\\
&\le C\,\lambda^\frac{p-2}{p}\left(\|\nabla\rho\|_{L^q}\|\nabla u\|_{L^\frac{pq}{q-p}}+\|\rho\dot u-b\cdot\nabla b\|_{L^p}\right),
\end{split}
\end{equation}
if $\|\mu(\rho_0)-1\|_{L^\infty}\le \varepsilon_4\le  \varepsilon_{41}$ for some small constant $\varepsilon_{41}$ owing to Lemma \ref{pl3}. Similarly, 
\begin{equation}
\begin{split}\label{bu1}
|\mathbf{I}_{12}(x)| 
&\leq C\,\|\mu(\rho)-1\|_{L^\infty}\int_\R\left(|\nabla\varphi_\lambda||x-y|^{-1}+(1-\varphi_\lambda)|x-y|^{-2}\right)|D u|(y)\,dy\\
& \leq C_3\,\|\mu(\rho_0)-1\|_{L^\infty}\left(\|\nabla u\|_{L^\infty}+\lambda^{-\frac{2(q-p)}{pq}}\|\nabla u\|_{L^\frac{pq}{q-p}}\right)\\
&\leq \frac{1}{2}\|\nabla u\|_{L^\infty}+C\,\varepsilon_4\,\lambda^{-\frac{2(q-p)}{pq}}\|\nabla u\|_{L^\frac{pq}{q-p}},
\end{split}
\end{equation}
provided that $ \varepsilon_4\le \varepsilon_{42} \triangleq (2C_3)^{-1}$. 
 
To bound $\mathbf{I}_2(x)=\nabla\Delta^{-1}\mathbb{P}(\rho\dot u-b\cdot\nabla b)$, we use  the Sobolev inequality to derive
\begin{equation}\label{C419}
\begin{split}
|\mathbf{I}_2(x)|&\leq C\left(\|\mathbf{I}_2\|_{L^4}+\|\nabla{\mathbf{I}}_2\|_{L^p}\right)\\
&\le C\left(\big\|\nabla^2\Delta^{-1}\mathbb{P}(\rho\dot u-b\cdot\nabla b)\big\|_{L^{\frac{4}{3}}}+\big\|\nabla^2\Delta^{-1}\mathbb{P}(\rho\dot u-b\cdot\nabla b)\big\|_{L^p}\right)\\
&\leq C\left(\|\rho\dot u-b\cdot\nabla b\|_{L^{\frac{4}{3}}}+\|\rho\dot u-b\cdot\nabla b\|_{L^p}\right),
\end{split}
\end{equation}
since the operator ``$\nabla^2\Delta^{-1}\mathbb{P}$" maps  $L^r$ ($1<r<\infty$) to itself (see \cite[Chapter 1]{S}). 

Substituting \eqref{Q1}--\eqref{C419} into \eqref{C422}  thus yields that if $\varepsilon_4\le\min\{\varepsilon_3,\varepsilon_{41},\varepsilon_{42}\}$,
\begin{equation}\label{C423}
\begin{split}
\|\nabla u\|_{L^\infty}
&\le  C\,\lambda^\frac{p-2}{p}\left(\|\nabla\rho\|_{L^q}\|\nabla u\|_{L^\frac{pq}{q-p}}+\|\rho\dot u-b\cdot\nabla b\|_{L^p}\right)\\
&\quad+C\,\varepsilon_4\,\lambda^{-\frac{2(q-p)}{pq}}\|\nabla u\|_{L^\frac{pq}{q-p}}+C\left(\|\rho\dot u-b\cdot\nabla b\|_{L^{\frac{4}{3}}}+\|\rho\dot u-b\cdot\nabla b\|_{L^p}\right)\\
&\leq C\left(\varepsilon_4^\frac{q(p-2)}{p(q-2)}\,\|\nabla\rho\|_{L^q}^\frac{2(q-p)}{p(q-2)}+1\right)\|\nabla u\|_{L^\frac{pq}{q-p}}\\
&\quad+C\left(\|\sqrt\rho\dot u\|_{L^2}+\|\nabla b\|_{L^4}+\|\rho\dot u\|_{L^{p}}\right),
\end{split}
\end{equation}
where we have chosen $\lambda=\varepsilon_4^{q/(q-2)}\left(1+\|\nabla\rho\|_{L^q}\right)^{-{q/(q-2)}}$ and  used \eqref{3.7}, \eqref{3.8}, \eqref{3.9}. 

Therefore, we deduce from \eqref{3.39} and \eqref{C423}  that 
\begin{equation}\begin{aligned}\label{C424}
\frac{d}{dt}\|\nabla\rho\|_{L^q}
&\leq C\,\|\nabla\rho\|_{L^q}\left(\varepsilon_4^\frac{q(p-2)}{p(q-2)}\|\nabla\rho\|_{L^q}^{\frac{2(q-p)}{p(q-2)}}+1\right)\\
&\quad\cdot\left(\|\nabla u\|_{L^\frac{pq}{q-p}}+\|\sqrt\rho\dot u\|_{L^2}+\|\nabla b\|_{L^4}+\|\rho\dot u\|_{L^{p}}\right).
\end{aligned}\end{equation}

By virtue of \eqref{3.7}--\eqref{3.8}, \eqref{q166}, \eqref{3.9}, \eqref{3.11}, and \eqref{p25}, we obtain that if $\varepsilon_4\le \varepsilon_{43}$ for some small constant $\varepsilon_{43}$ owing to Lemma \ref{pl3},
\begin{equation}\begin{aligned}\label{3.22}
&\int_0^T\|\nabla u\|_{L^\frac{pq}{q-p}}\,dt\\
&\le C\int_0^T\|\nabla u\|_{L^2}^{\frac{2(q-p)}{pq}}\|\rho\dot u-b\cdot\nabla b\|_{L^2}^{1-\frac{2(q-p)}{pq}}\,dt\\
&\leq C\int_0^{\sigma(T)}\left(\|\nabla u\|_{L^2}+\|\sqrt\rho\dot u\|_{L^2}+\||\nabla b|\,|b|\|_{L^2}\right)\,dt\\
&\quad+C\int_{\sigma(T)}^T\left(t^{\eta_1}\|\nabla u\|_{L^2}^2\right)^{\frac{q-p}{pq}}\left(t^{1+\eta_1}\|\rho\dot u-b\cdot\nabla b\|_{L^2}^2\right)^{\frac{1}{2}-\frac{q-p}{pq}}t^{-\frac{1+\eta_1}{2}+\frac{q-p}{pq}}\,dt\\
&\leq C+C\int_{\sigma(T)}^Tt^{-(1+\eta_1)+\frac{2(q-p)}{pq}}\,dt\le C,
\end{aligned}\end{equation}
since $1+\eta_1-{2(q-p)}/{(pq)}>1$ by taking $\eta_1=1/2+1/p-1/q\in(1/2,1)$. 

Then it follows from \eqref{2.1} with $\xi=1/24$, \eqref{3.7}, \eqref{k37}, and \eqref{k38} that for any $2\le r<72/5$,
\begin{equation}\begin{aligned}\label{3.23}\notag
\|\rho\dot u\|_{L^r} 
&\le C\,(1+t)^{1/12}\left(N(1+t)\right)^2 \left(\|\sqrt\rho \dot u\|_{L^2}+\|\nabla\dot u\|_{L^2}\right)\\
&\le C\,(1+t)^{25/12}\left(\|\sqrt\rho \dot u\|_{L^2}+\|\nabla\dot u\|_{L^2}\right),
\end{aligned}\end{equation}
which together with \eqref{3.7}, \eqref{3.9}, \eqref{3.10}, and \eqref{3.11} gives that 
\begin{equation}\begin{aligned}\label{3.24}
&\int_0^T\|\rho\dot u\|_{L^p}\,dt\\
&\le\int_0^T\|\rho\dot u\|_{L^{2}}^\frac{12-p}{5p}\|\rho\dot u\|_{L^{12}}^\frac{6(p-2)}{5p}\,dt\\
&\le C\int_0^T (1+t)^\frac{5(p-2)}{2p}\|\sqrt\rho\dot u\|_{L^2}^\frac{12-p}{5p}\left(\|\sqrt\rho\dot u\|_{L^2}+\|\nabla\dot u\|_{L^2}\right)^\frac{6(p-2)}{5p}\,dt\\
&\le C\int_0^{\sigma(T)}\left(\|\sqrt\rho\dot u\|_{L^2}+\|\sqrt\rho\dot u\|_{L^2}^\frac{12-p}{5p}\left(t\|\nabla\dot u\|_{L^2}^2\right)^\frac{3(p-2)}{5p}t^{-\frac{3(p-2)}{5p}}\right)\,dt\\
&\quad+ C\int_{\sigma(T)}^Tt^{\frac{5(p-2)}{2p}-\frac{1+\eta_2}{2}}\left(t^{1+\eta_2}\|\sqrt\rho\dot u\|_{L^2}^2+t^{2+\eta_2}\|\nabla\dot u\|_{L^2}^2\right)^\frac{1}{2}\,dt\\
&\leq C+C\int_0^{\sigma(T)}t^{-\frac{6(p-2)}{5p}}\,dt+ C\int_{\sigma(T)}^Tt^{-(1+\eta_2)+\frac{5(p-2)}{p}}\,dt  \leq C,
\end{aligned}\end{equation}
since $6(p-2)/5p<1$ and $1+\eta_2-5(p-2)/p>1$ by taking $\eta_2=2-3/p\in(1/2,1)$.

Furthermore, we obtain from \eqref{3.8},  \eqref{3.9},    and \eqref{3.11} with $\eta=5/6$ that
\begin{equation}\begin{aligned}\label{3.25}
&\int_0^{T}\left(\|\sqrt\rho\dot u\|_{L^2}+\|\nabla b\|_{L^4}\right)\,dt\\
&\le C\int_0^{\sigma(T)}\left(\|\sqrt\rho\dot u\|_{L^2}+\|\nabla b\|_{L^2}+\|\Delta b\|_{L^2}\right)\,dt\\
&\quad +C\int^T_{\sigma(T)}\left(t^{-\frac{11}{12}}\left(t^\frac{11}{6}\|\sqrt\rho\dot u\|_{L^2}^2\right)^\frac{1}{2}+t^{-\frac{2}{3}}\left(t^\frac{5}{6}\|\nabla b\|_{L^2}^2\right)^\frac{1}{4}\left(t^\frac{11}{6}\|\Delta b\|_{L^2}^2\right)^\frac{1}{4}\right)\,dt\\
&\leq C+C\int^T_{\sigma(T)}\left(t^{-\frac{11}{6}}+t^{-\frac{4}{3}}\right)\,dt\le C.
\end{aligned}\end{equation}

As a result, applying Gronwall's inequality to \eqref{C424}  and using \eqref{3.22}--\eqref{3.25} infer that
$$\sup_{ t\in[0,T]}\|\nabla\rho\|_{L^q}\leq  \text{exp}\left\{C_4\left(\varepsilon_4^\frac{q(p-2)}{p(q-2)}\sup_{ t\in[0,T]}\|\nabla\rho\|_{L^q}^{\frac{2(q-p)}{p(q-2)}}+1\right)\right\}\|\nabla\rho_0\|_{L^q},$$
which,  by a standard continuity argument, implies 
\begin{equation}\label{nr}
\sup_{ t\in[0,T]}\|\nabla\rho\|_{L^q}\leq 3e^{C_4}\|\nabla\rho_0\|_{L^q},
\end{equation}
by supposing $\|\mu(\rho_0)-1\|_{L^\infty}\leq\varepsilon_{4}\triangleq \min\{(6\|\nabla\rho_0\|_{L^q})^{-\frac{2(q-p)}{q(p-2)}}e^{-C_4\frac{p+2}{p-2}},\,\varepsilon_3,\,\varepsilon_{41},\,\varepsilon_{42},\,\varepsilon_{43}\}$. 

Finally, in the similar way based on \eqref{C423} and  \eqref{3.22}--\eqref{nr}, we conclude that
\begin{equation}\notag
\sup_{ t\in[0,T]}\|\nabla \rho\|_{L^2}+\int_0^T\|\nabla u\|_{L^\infty}\,dt\leq C.
\end{equation}
This combined with \eqref{3.7} and \eqref{nr} gives \eqref{3.19} and therefore finishes the proof of Lemma \ref{le6}.
\end{proof}

\begin{lemma}\label{le7}
Under the conditions of Theorem \ref{T1}, suppose that
$$\|\mu(\rho_0)-1\|_{L^\infty}\leq\varepsilon_0,$$
with $\varepsilon_0$ $(\le \varepsilon_4)$  depending on  $\bar\rho$, $\mu$, $q$, and the initial values. For any $\eta\in(0,1)$, there exists a positive constant $C$ depending on  $\eta$, $\bar\rho$, $\mu$, $q$, 
and the initial values such that 
\begin{equation}\begin{aligned}\label{3.288}
&\sup_{t\in[0,T]}t^{1+\eta}\left(\|\nabla u\|_{H^1}^2+\| P\|_{H^1}^2+\|\nabla b\|_{H^1}^2+\|\sqrt\rho  u_t\|_{L^2}^2+\|b_t\|_{L^2}^2\right)\\
&+\sup_{t\in[0,T]}\|\rho\,\tilde x^2\|_{L^1}+\int_0^T(1+t^\eta)\left(\|\nabla u\|_{H^1}^2+\|P\|_{H^1}^2+\|\sqrt\rho  u_t\|_{L^2}^2+\|b_t\|_{L^2}^2\right)\,dt\\
&+\int_0^Ta(t)\left(\|\nabla  u\|_{W^{1,q}}^2+\| P\|_{W^{1,q}}^2\right)\,dt\\
&+\int_0^Tt^{1+\eta}\left(\|\nabla u_t\|_{L^2}^2+\|\nabla b_t\|_{L^2}^2+\|\nabla^3 b\|_{L^2}^2\right)\,dt\leq C,
\end{aligned}\end{equation}
where $a(t)\triangleq\min\{t,t^\eta\}$, and for any $\alpha\in(0,1/2)$,  there exists a positive constant $C$ depending on $\alpha$, $\eta$, $\bar\rho$, $\mu$, 
and the initial values such that for a.e. $t>0$,
\begin{equation}\label{uui}
\| u-u_\infty\|_{L^{2/\alpha}}+|u_\infty(t)|\le C\,t^{-(1+\eta)/2}.
\end{equation}
\end{lemma}
\begin{proof}
In fact, we deduce from the proof of Lemma \ref{le11} that
\begin{equation}\label{u12}\begin{aligned}
\|u-u_\infty\|_{L^{8}}&\le C\left(1+\|\rho\,\tilde x^2\|_{L^1}\right)^{1/6}\,\|\sqrt\rho\dot u\|_{L^2}+C\left(\|b\|_{L^{16/5}}^2+\|b\|_{L^4}^2\right)\\
&\le C\left(1+\|\rho\,\tilde x^2\|_{L^1}\right)^{1/6}\,\|\sqrt\rho\dot u\|_{L^2}\\
&\quad+C\left(\|b\|_{L^{16/15}}^{2/3}\|\nabla b\|_{L^2}^{4/3}+\|b\|_{L^2}\|\nabla b\|_{L^2}\right),
\end{aligned}\end{equation}
which along with \eqref{3.7},  \eqref{3.9}, \eqref{p22},  \eqref{k37}, and \eqref{3.11} with $\eta=5/6$ leads to 
\begin{equation}\notag\begin{aligned}
\int_0^T\|\sqrt\rho u\|_{L^2}\,dt
&\le \int_0^T\|\rho\|_{L^{4/3}}^{1/2}\|u-u_\infty\|_{L^{8}}\,dt\\
&\le C\int_0^T\left((1+t)^{1/3}\,\|\sqrt\rho\dot u\|_{L^2}+\|b\|_{L^{16/15}}^{2/3}\|\nabla b\|_{L^2}^{4/3}+\|b\|_{L^2}\|\nabla b\|_{L^2}\right)\,dt\\
&\le C\int_0^{\sigma(T)}\left(\|\sqrt\rho\dot u\|_{L^2}+\|\nabla b\|_{L^2}^{4/3}+\|\nabla b\|_{L^2}\right)\,dt\\
&\quad+C\int_{\sigma(T)}^T\left(t^{1/3-17/12}+t^{-11/9}+t^{-4/3}\right)\,dt\le C.
\end{aligned}\end{equation}
Therefore, by the proof of Lemma \ref{le9}, we have
\begin{equation}\label{C426}
\sup_{t\in[0,T]}\|\rho\,\tilde x^2\|_{L^1}^{1/2}\leq \|\rho_0\,\tilde x^2\|_{L^1}^{1/2}+2\int_0^T\|\sqrt\rho u\|_{L^2}\,dt\leq C,
\end{equation}
and there is a constant $\tilde N>0$ independent of $t$ such that 
\begin{equation}\notag
\int_{B_{\tilde N}}\rho(x,t)\,dx\geq\frac{1}{2}\|\rho_0\|_{L^1},\quad t\in(0,T].
\end{equation}
Based on these,  the proof of Lemma \ref{jq} with a little modification implies that for any function $f$ satisfying $\sqrt\rho f,\nabla f\in L^2$ and $p\in [2,\infty)$,
\begin{equation}\label{un}
\|\rho f\|_{L^p}\le C\left(\|\sqrt\rho f\|_{L^2}+\|\nabla f\|_{L^2}\right).
\end{equation}

Next, it follows from  \eqref{C225}, \eqref{C202}, \eqref{1.5}, \eqref{3.8}, \eqref{har}, and \eqref{C426} that
\begin{equation}\label{pe}\begin{aligned}
\|P\|_{L^2}
&\le \|\Delta^{-1}\mathrm{div}(\rho\dot u-b\cdot\nabla b)\|_{L^2}+\|\Delta^{-1}\mathrm{div}\mathrm{div}(2\mu(\rho)Du)\|_{L^2}\\
&\le C\,\|\nabla\Delta^{-1}\mathrm{div}(\rho\dot u)\|_{L^1}+ \|\Delta^{-1}\mathrm{div}\mathrm{div}(b\otimes b)\|_{L^2}+C\,\|\nabla u\|_{L^2}\\
&\le C\left(\|\rho\dot u\|_{\mathcal{H}^1}+\|b\|_{L^4}^2+\|\nabla u\|_{L^2}\right)\\
&\le C\left(\|\sqrt\rho\dot u\|_{L^2}+\|\nabla b\|_{L^2}+\|\nabla u\|_{L^2}\right).
\end{aligned}\end{equation}
This along with Lemma \ref{pl3}, \eqref{3.7}, \eqref{1.5},  \eqref{pe}, \eqref{3.19}, \eqref{q166}, \eqref{3.9}, \eqref{C423}, and \eqref{un} shows that there is a  positive constant $\varepsilon_{0}\,(\le \varepsilon_4)$ such that if $\|\mu(\rho_0)-1\|_{L^\infty}\leq\varepsilon_{0},$ 
$$\|\nabla u\|_{H^1}+\|P\|_{H^1}\le C\left(\|\sqrt\rho\dot u\|_{L^2}+\|\nabla b\|_{L^2}+\||\nabla b||b|\|_{L^2}+\|\nabla u\|_{L^2}\right),$$
$$\|\nabla u\|_{W^{1,q}}+\|P\|_{W^{1,q}}\le C\left(\|\sqrt\rho \dot u\|_{L^2}+\|\nabla\dot u\|_{L^2}+\|\nabla b\|_{L^2}+\|\Delta b\|_{L^2}+\|\nabla u\|_{L^2}\right).$$
Then we deduce from the above estimates and Lemmas \ref{le2}--\ref{le4}, \ref{le5} that for $a(t)\triangleq\min\{t,t^\eta\}$,
\begin{equation}\begin{aligned}\label{3.40}
&\sup_{t\in[0,T]}t^{1+\eta}\left(\|\nabla u\|_{H^1}^2+\|P\|_{H^1}^2\right)+\int_0^T (1+ t^{\eta})\left(\|\nabla u\|_{H^1}^2+\|P\|_{H^1}^2\right)\,dt\\
&+\int_0^Ta(t)\left(\|\nabla u\|_{W^{1,q}}^2+\|P\|_{W^{1,q}}^2\right)\,dt\leq C.
\end{aligned}\end{equation}

Moreover, notice that \eqref{u12}, \eqref{C426}, \eqref{p22}, and \eqref{3.9} indicate
\begin{equation}\begin{aligned}\label{uin}
\|u\|_{L^\infty}
&\le \|u-u_\infty\|_{L^{\infty}}+|u_\infty|\\
&\le C\,\|u-u_\infty\|_{L^{8}}^{4/5}\|\nabla^2 u\|_{L^2}^{1/5}+C\left(\|\sqrt\rho\dot u\|_{L^2}+\|\nabla b\|_{L^2}\right)\\
&\le C\left(\|\sqrt\rho\dot u\|_{L^2}+\|\nabla b\|_{L^2}+\|\nabla^2 u\|_{L^2}\right).
\end{aligned}\end{equation}
Direct calculation along with  \eqref{3.7}, \eqref{uin},  and \eqref{3.9} yields
\begin{equation*}\begin{aligned}
\|\sqrt\rho u_t\|_{L^2}
&\le \|\sqrt\rho \dot u\|_{L^2}+\|\sqrt\rho u\cdot\nabla u\|_{L^2}\\
&\le \|\sqrt\rho \dot u\|_{L^2}+\|\rho\|_{L^{\infty}}^{1/2}\|u\|_{L^{\infty}}\|\nabla u\|_{L^{2}}\\
&\le C\left(\|\sqrt\rho\dot u\|_{L^2}+\|\nabla b\|_{L^2}+\|\nabla^2 u\|_{L^2}\right).
\end{aligned}\end{equation*}
Similarly, 
\begin{equation*}\begin{aligned}
\|\nabla u_t\|_{L^2}
&\le \|\nabla \dot u\|_{L^2}+\|\nabla u\|_{L^4}^2+\||u|\,|\nabla^2 u|\|_{L^2}\\
&\le \|\nabla \dot u\|_{L^2}+C\,\|\nabla u\|_{L^{2}}\|\nabla^2 u\|_{L^{2}}+C\,\|u\|_{L^\infty}\|\nabla^2 u\|_{L^2}\\
&\le \|\nabla \dot u\|_{L^2}+C\left(\|\nabla u\|_{H^{1}}^2+\|\sqrt\rho \dot u\|_{L^2}^2+\|\nabla b\|_{L^2}^2\right).
\end{aligned}\end{equation*}
By using these estimates, \eqref{3.40}, and Lemmas \ref{le2}--\ref{le4}, \ref{le5}, we get that
\begin{equation}\begin{aligned}\label{3.43}
\sup_{t\in[0,T]}t^{1+\eta}\,\|\sqrt\rho  u_t\|_{L^2}^2+\int_0^T\left((1+t^{\eta})\|\sqrt\rho u_t\|_{L^2}^2+t^{1+\eta}\|\nabla u_t\|_{L^2}^2\right)\,dt\le C.
\end{aligned}\end{equation}

Furthermore, integrating \eqref{C101}$_3$ multiplied by $2\Delta b_t$ by parts gives that 
\begin{equation}\begin{aligned}\label{bt}\notag
&\frac{d}{dt}\|\Delta b\|_{L^2}^2+\|\nabla b_t\|_{L^2}^2+C\|\nabla^3 b\|_{L^2}^2\\
&\le C\left(\|\nabla u\|_{L^4}^2\|\nabla b\|_{L^4}^2+\|u\|_{L^\infty}^2\|\nabla^2 b\|_{L^2}^2+\|b\|_{L^\infty}^2\|\nabla^2 u\|_{L^2}^2\right)\\
&\le C\,\|\nabla u\|_{L^2}\|\nabla^2 u\|_{L^2}\|\nabla b\|_{L^2}\|\nabla^2 b\|_{L^2}+C\,\|u\|_{L^\infty}^2\|\nabla^2 b\|_{L^2}^2\\
&\quad+C\,\|b\|_{L^2}\|\nabla^2 b\|_{L^2}\|\nabla^2 u\|_{L^2}^2.
\end{aligned}\end{equation}
Multiplying it by $t^{1+\eta}$ and using \eqref{uin},  \eqref{3.40}, Lemmas \ref{le2}--\ref{le4}, \ref{le5}  lead to
\begin{equation}\begin{aligned}\label{3.44}
\sup_{t\in[0,T]}t^{1+\eta}\,\|\Delta b\|_{L^2}^2+\int_0^Tt^{1+\eta}\left(\|\nabla b_t\|_{L^2}^2+\|\nabla^3 b\|_{L^2}^2\right)\,dt\le C.
\end{aligned}\end{equation}
In addition, observe that
\begin{equation}\begin{aligned}\notag
\|b_t\|_{L^2}
&\le C\left(\|\Delta b\|_{L^2}+\|u\cdot\nabla b\|_{L^2}+\|b\cdot\nabla u\|_{L^2}\right)\\
&\le C\left(\|\Delta b\|_{L^2}+\|u\|_{L^\infty}\|\nabla b\|_{L^2}+\|b\|_{L^\infty}\|\nabla u\|_{L^2}\right)\\
&\le C\left(\|\Delta b\|_{L^2}+\|\sqrt\rho\dot u\|_{L^2}+\|\nabla b\|_{L^2}+\|\nabla u\|_{L^2}\right),
\end{aligned}\end{equation}
which combined with Lemmas \ref{le2}--\ref{le4}, \ref{le5}, and \eqref{3.40}--\eqref{3.44} concludes \eqref{3.288}. Finally, for any $\alpha\in(0,1/2)$,  it follows from \eqref{p10}, \eqref{C426}, and \eqref{3.11} that  \eqref{uui} holds. 
The proof of Lemma \ref{le7} is finished.

\end{proof}

\section{Proof of main Theorem}\label{S5}
\quad With all a priori estimates in Section \ref{S3} at hand, we can now prove the global existence result stated in Theorem \ref{T1} by the standard extension method.

{\it Proof of Theorem \ref{T1}:} 
Besides the conditions of Theorem \ref{T1}, we assume further that  \eqref{com} holds, that is 
\begin{equation*}
\|\mu(\rho_0)-1\|_{L^\infty}\leq\varepsilon_0,
\end{equation*}
with $\varepsilon_0$ as in Lemma \ref{le7}. The local existence and uniqueness of the strong solutions are guaranteed by Lemma \ref{L41}.  Let $T^*>0$ be the maximal  existence time of the strong solution $(\rho,u,P,b)$. Our aim is to prove $T^\ast=\infty$. To do this, we use contradiction arguments. More precisely, we assume  $0<T^*<\infty$. Then based on \eqref{1.7}, we obtain that for any $0< \tau<T< T^{*}$,
\begin{equation}\label{C502}
\rho \in C([0,T];L^1\cap H^1\cap W^{1,q}),
\end{equation}
\begin{equation}\label{C501}
\rho u,\,\nabla u,b,\nabla b\in L^{\infty}(\tau,T;H^1)\cap H^1(\tau,T;L^2) \subset C([\tau,T];L^2\cap L^q).
\end{equation} 
With the uniform (with respect to $T$) estimates in Lemmas \ref{le1}, \ref{le5}--\ref{le7} and  \eqref{C501}--\eqref{C502}, we can define new initial data at $T^\ast$,
$$(\rho, u,b)(x,T^*)=\lim_{t\rightarrow T^*}(\rho, u,b)(x,t).$$
The limiting data satisfy all the assumptions of Lemma \ref{L41}.  
Hence the local solution can be continued beyond \(T^*\), which is the desired contradiction.  
Finally,  the additional properties of solutions \eqref{C102}--\eqref{1.9} follow from Lemmas \ref{le5}--\ref{le7} directly. The proof of Theorem \ref{T1} is completed.  
\thatsall

\section*{Statements and Declarations}
{\bf Conflict of interest statement}
On behalf of all authors, the corresponding author states that there is no conflict of interest. \\
{\bf  Data availability statement} 
No data were used for this research.


\end{document}